\RequirePackage{fix-cm}
\documentclass[10pt]{article}

\usepackage{latex_styles_common}
\usepackage{cahier}

\newcommand{\cahiernumber}{53}  

\def\papertitle{An Exact Penalty Method for Large-Scale Equality-Constrained Nonlinear Optimization}
\def\authorone{Maxence Gollier}
\def\authortwo{Youssef Diouane}
\def\authorthree{Dominique Orban}

\hypersetup{
  colorlinks=true, linkcolor=blue, citecolor=blue, urlcolor=blue,
  pdftitle={An Exact Penalty Method for Large-Scale Equality-Constrained Nonlinear Optimization},
  pdfauthor={\authorone, \authortwo and \authorthree},
  pdfsubject={Large-Scale Exact Penalty Methods},
  pdfkeywords={Nonlinear Optimization, Nonconvex optimization, Constrained Optimization, Complexity Analysis, Julia (Programming Language)},
}

\title{\papertitle}

\author{
  \authorone\thanks{%
    GERAD and Department of Mathematics and Industrial Engineering, Polytechnique Montr\'eal.
    E-mail: \mailto{maxence-2.gollier@polymtl.ca}.
    Research partially supported by a PBEEE Scholarship.}
  \and
  \authortwo\thanks{%
    GERAD and Department of Mathematics and Industrial Engineering, Polytechnique Montr\'eal.
    E-mail: \mailto{youssef.diouane@polymtl.ca}.
    Research partially supported by an NSERC Discovery Grant.}
  \and
  \authorthree\thanks{%
    GERAD and Department of Mathematics and Industrial Engineering, Polytechnique Montr\'eal.
    E-mail: \mailto{dominique.orban@gerad.ca}.
    Research partially supported by an NSERC Discovery Grant.}
}
\date{\today}

\begin{document}

\maketitle
\thispagestyle{cahier}

\begin{abstract}
    We present an exact penalty method for equality-constrained nonlinear optimization.
    Our approach transforms the smooth, constrained problem into a sequence of nonsmooth, unconstrained problems by penalizing the constraint violation by its Euclidean norm.
    An outer loop updates the penalty parameter and the subproblem tolerances, while an inner loop solves each penalized subproblem with a quadratic regularization method.
    The method also includes a local infeasibility detection mechanism, an extrapolation strategy for warm-starting subproblems, and nonmonotone acceptance rules.
    We establish global convergence to first-order stationary points under a bounded penalty parameter assumption, with a worst-case complexity of \(\mathcal{O}(|\log \epsilon|\,\epsilon^{-2})\) inner iterations.
    We implement the method in Penelopt, an open-source Julia solver for constrained nonlinear optimization.
    On \(224\) problems from the CUTEst collection, Penelopt is competitive with IPOPT and more robust than the recent solvers MadNLP and Uno.
    These results establish exact penalty methods as a practical alternative to conventional methods for solving equality-constrained nonlinear optimization problems.
    To the best of our knowledge, Penelopt is the first practical implementation of an exact penalty method for large-scale constrained nonlinear optimization.
\end{abstract}




\section{Introduction}%
\label{sec:introduction}

We present a nonsmooth exact penalty method to solve the general equality-constrained nonlinear optimization problem
\begin{equation}%
  \label{eq:nlp}
  \minimize{x \in \R^n} \ f(x) \quad \st \ c(x) = 0,
\end{equation}
where \(f: \R^n \to \R\) and \(c: \R^n \to \R^m\), with \(m \leq n\), are continuously differentiable and may be nonconvex.

Following the approach of the penalty methods presented in \citep{diouane-gollier-orban-2026}, our algorithm computes a sequence of approximate solutions of nonsmooth penalized problems
\begin{equation}%
  \label{eq:penalty-nlp}
  \minimize{x \in \R^n} \ f(x) + \tau \|c(x)\|_2,
\end{equation}
for a succession of non-decreasing penalty parameters \(\tau > 0\).

Our work builds directly upon that of \citet{diouane-gollier-orban-2026}, who propose the first practical implementation of the \(\ell_2\)-penalty method for solving small-scale instances of~\eqref{eq:nlp}, relying on the proximal methods of \citet{aravkin-baraldi-orban-2022} and \citet{diouane-habiboullah-orban-2026} to handle the resulting nonsmooth penalized subproblems~\eqref{eq:penalty-nlp}.
The present work significantly extends both the algorithmic design and the practical implementation of \citet{diouane-gollier-orban-2026}.
We emphasize the following contributions:
\begin{enumerate}%
  \item We introduce a new step computation and acceptance rule that substantially reduces the number of linear system solves required by the algorithm compared to \citep{diouane-gollier-orban-2026};
  \item We incorporate several practical heuristics that improve performance on a broad range of problems;
  \item We provide a thoroughly tested and well-documented open-source Julia \citep{bezanson-edelman-karpinski-shah-2017} implementation named \emph{Penelopt}.
\end{enumerate}
Together, these advances allow us to draw significantly stronger conclusions than those of \citet{diouane-gollier-orban-2026}.
In particular, our numerical results show that Penelopt is competitive with the state-of-the-art solver IPOPT on large-scale problems, establishing exact \(\ell_2\)-norm penalty methods as a practical alternative to SQP methods for solving~\eqref{eq:nlp}.
Moreover, to the best of our knowledge, Penelopt is the first practical implementation of an exact penalty method for large-scale constrained nonlinear optimization.

The paper is organized as follows.
In \Cref{sec:background}, we review key concepts and results that underpin our algorithm.
In \Cref{sec:algorithm-description}, we describe the overall algorithm.
The convergence properties of the method are then established in \Cref{sec:convergence-analysis} under mild assumptions, followed by the presentation of the step computation procedure in \Cref{sec:step-computation}.
In \Cref{sec:implementation-details}, we discuss the practical aspects of the implementation, including nonmonotone strategies and other heuristics designed to improve performance.
\Cref{sec:software} describes the design and testing of our software.
Finally, in \Cref{sec:numerical-results}, we evaluate the numerical efficiency and robustness of our solver on CUTEst instances \citep{gould-orban-toint-2015}, and compare its performance with the state-of-the-art solver IPOPT \citep{wachter-biegler-2006} as well as the recent Uno \citep{vanaret-leyffer-2026, vanaret-montoison-2026} and the pure Julia MadNLP \citep{shin-anitescu-pacaud-2024} solvers.

\subsection*{Related Work}%
\label{sec:related-work}

The pioneering work of \citet{pietrzykowski-1969} was the first to establish the theoretical advantages of nonsmooth penalty functions, with particular emphasis on the \(\ell_1\)-norm penalty.
His foundational results were subsequently extended and strengthened by \citet{charalambous-1978}, \citet{han-mangasarian-1979}, \citet{coleman-conn-1980}, and \citet{bazaraa-goode-1982}.

More recently, \citet{cartis-gould-toint-2011} introduced trust-region and quadratic regularization algorithms for solving general \emph{composite} optimization problems of the form
\begin{equation}
  \label{eq:composite-objective}
  \minimize{x \in \R^n} f(x) + h(c(x)),
\end{equation}
where \(h\) is convex and Lipschitz continuous, encompassing nonsmooth penalty formulations as a special case.
Unlike the present work, \citet{cartis-gould-toint-2011} do not provide a practical implementation of their methods.
Independently, \citet{aravkin-baraldi-orban-2022} and \citet{diouane-habiboullah-orban-2026} developed globally convergent proximal methods capable of solving~\eqref{eq:composite-objective}, in which steps are computed by considering a model of \(h \circ c\).
Then, \citet{diouane-gollier-orban-2026} proposed an \(\ell_2\)-norm penalty approach for~\eqref{eq:nlp}, in which the resulting penalized subproblems are handled by such proximal methods, as discussed above.
In addition, \citet{diouane-gollier-orban-2026} provide a first practical implementation for small-scale problems and describe, contrary to \citep{cartis-gould-toint-2011}, a practical step computation procedure to solve~\eqref{eq:penalty-nlp}.
We now briefly review alternative state-of-the-art methods for solving~\eqref{eq:nlp}.

The augmented Lagrangian, originally due to \citet{hestenes-1969} and \citet{powell-1969}, is a prominent framework for constrained optimization.
It applies a quadratic penalty to the Lagrangian of~\eqref{eq:nlp}.
Its theory is developed in depth by \citet{bertsekas-1982}.
State-of-the-art implementations include ALGENCAN \citep{andreani-birgin-martinez-schuverdt-2008,birgin-martinez-2014}, LANCELOT \citep{conn-gould-toint-1991,conn-gould-toint-1992}, and MINOS \citep{murtagh-saunders-1978,murtagh-saunders-1982}.

Sequential quadratic programming (SQP) constitutes another major paradigm for~\eqref{eq:nlp}.
It solves a sequence of quadratic subproblems that approximate~\eqref{eq:nlp}.
The SQP framework was first introduced by \citet{wilson-1963} and subsequently popularized by \citet{garcia-palomares-mangasarian-1976} and \citet{han-1977}.
An overview of SQP is provided by \citet{boggs-tolle-1995}.
Widely used implementations include SNOPT \citep{gill-murray-saunders-2002}, FilterSQP \citep{fletcher-leyffer-1998}, and, more recently, Uno \citep{vanaret-leyffer-2026,vanaret-montoison-2026}.
In addition, the IPOPT \citep{wachter-biegler-2006}, KNITRO \citep{byrd-hribar-nocedal-1999,waltz-morales-nocedal-orban-2006}, and MadNLP \citep{shin-anitescu-pacaud-2024} solvers reduce to SQP methods when applied to~\eqref{eq:nlp}.

\subsection*{Notation}%
\label{sec:notation}
The identity matrix of size \(n\) is denoted by \(I_n\), or \(I\) if the context is clear.
We use \(J(x)\) to denote the Jacobian of \(c\) at \(x\).
We call a function \(\mathcal{C}^1\) if it is continuously differentiable.
For \(p \geq 1\), the vector \(\ell_p\)-norm is \(\|\cdot\|_p\).
For matrices, \(\|\cdot\|_p\) denotes the operator norm induced by the vector \(\ell_p\)-norm.
We use \(\lambda_{\min}(A)\) and \(\lambda_{\max}(A)\) to denote the smallest and largest eigenvalues of a symmetric matrix \(A\), respectively.
We write \(A \succ 0\) if \(A\) is positive definite, i.e., \(\lambda_{\min}(A) > 0\), and \(A \succeq 0\) if \(A\) is positive semidefinite, i.e., \(\lambda_{\min}(A) \geq 0\).
If \(g\), \(h\) are two positive functions of \(\epsilon > 0\), the notation \(g = \mathcal{O}(h)\) means that there is \(C > 0\) such that \(\limsup_{\epsilon \to 0} g(\epsilon)/h(\epsilon) \leq C\).
When clear from the context, given \(x_k \in \R^n\), we write \(f_k\) for \(f(x_k)\), \(\nabla f_k\) for \(\nabla f(x_k)\), \(c_k\) for \(c(x_k)\) and \(J_k\) for \(J(x_k)\).
We refer to \(\epsilonm\) as the machine precision for the finite precision arithmetic.
All numerical results presented in this paper are performed in double precision with \(\epsilonm \approx 10^{-16}\).

\section{Background}%
\label{sec:background}

We say that the Mangasarian-Fromovitz constraint qualification (MFCQ) holds at \(\bar x \) for~\eqref{eq:nlp} if \(J(\bar x)\) has full row rank.

Recall that under the MFCQ, every local solution \(\bar x\) of~\eqref{eq:nlp} satisfies the first-order necessary optimality conditions
\begin{subequations}%
  \label{eq:kkt}
  \begin{align}%
    \label{eq:kkt-dual}
    \nabla f(\bar x) + { J(\bar x) }^T \bar y & = 0, \\
    \label{eq:kkt-primal}
    c(\bar x)                                 & = 0,
  \end{align}
\end{subequations}
where \(\bar y \in \R^m\) is the unique vector of Lagrange multipliers associated with \(\bar x\).

We say that \(\bar x \in \R^n\) is a strict local minimizer of~\eqref{eq:nlp} if \(c(\bar x) = 0\) and there exists an open set \(\mathcal{V}\) with \(\bar x \in \mathcal{V}\) such that \(f(\bar x) < f(x)\) for all \(x \in \mathcal{V} \setminus \{\bar x\}\) satisfying \(c(x) = 0\).
Exactness of the \(\ell_2\)-norm penalty function means the following.

\begin{theorem}[\protect{\citealp[Theorem 4.4]{han-mangasarian-1979}}]%
  \label{thm:exactness}
  Let \(\bar x\) be a strict local minimizer of~\eqref{eq:nlp} that satisfies the MFCQ\@.
  Then, there exists a finite threshold \(\bar \tau > 0\) such that \(\bar x\) is a local minimizer of~\eqref{eq:penalty-nlp} for all \(\bar \tau \leq \tau < \infty\).
  Moreover, \(\bar \tau\) satisfies \(\bar \tau \geq \|\bar y\|_2\), where \(\bar y\) is the unique vector of Lagrange multipliers associated with \(\bar x\).
\end{theorem}

\Cref{thm:exactness} motivates the use of the \(\ell_2\)-norm penalty in~\eqref{eq:penalty-nlp}.
Contrary to inexact penalty functions, such as the squared \(\ell_2\)-norm, we do not need to drive \(\tau \to \infty\) to establish convergence to a local solution of~\eqref{eq:nlp}.

We now consider the first-order necessary optimality conditions for~\eqref{eq:penalty-nlp}.
\begin{proposition}%
  \label{prop:first-order-penalized}
  Let \(\bar x\) be a local solution of~\eqref{eq:penalty-nlp} for some \(\tau > 0\).
  Then, there exists \(\bar y \in \R^m\) such that
  \begin{equation}%
    \label{eq:outer-kkt}
    \nabla f(\bar x) + { J(\bar x) }^T \bar y = 0.
  \end{equation}
  Moreover, there is \(\bar \alpha \geq 0\) such that
  \begin{equation}%
    \label{eq:inner-kkt}
    c(\bar x) - \bar \alpha \bar y = 0, \quad \bar \alpha(\|\bar y\|_2^2 - \tau^2) = 0, \quad \|\bar y\|_2 \leq \tau.
  \end{equation}
\end{proposition}
\begin{proof}%
  We rewrite~\eqref{eq:penalty-nlp} as
  \begin{equation}%
    \label{eq:min-max}
    \min_{x} f(x) + \tau \|c(x)\|_2 = \min_x \max_{\|y\|_2 \leq \tau} \ f(x) + y^T c(x).
  \end{equation}
  Convexity of the maximization problem in~\eqref{eq:min-max} implies that its necessary and sufficient optimality conditions are~\eqref{eq:inner-kkt}.
  Then, the necessary first-order optimality conditions for the minimization problem in~\eqref{eq:min-max} are~\eqref{eq:outer-kkt}.
\end{proof}

\section{Algorithm Description}%
\label{sec:algorithm-description}

\subsection{The Exact Penalty Method}%
\label{sec:exact-penalty-method}

Our algorithm iteratively computes an approximate solution of~\eqref{eq:penalty-nlp} for a fixed penalty parameter \(\tau_k\), then either increases \(\tau_k\) or improves the accuracy of the current solution.
We call this process the ``outer loop'' of our algorithm and use the index \(k\) to denote the iteration counter.

At each outer iteration \(k\), we compute primal and dual variables \((x_k, y_k)\) that approximately satisfy~\eqref{eq:outer-kkt}.
When the pair \((x_k, y_k)\) additionally (approximately) satisfies the first-order necessary optimality conditions of~\eqref{eq:nlp}, we declare first-order stationarity and terminate the algorithm.
In view of~\eqref{eq:kkt}, our algorithm stops when \((x_k, y_k)\) satisfies
\begin{equation}%
  \label{eq:overall-stopping-criterion}
  \| \nabla f(x_k) + { J(x_k) }^T y_k \|_\infty \leq \epsilon_D \quad \text{and} \quad \|c(x_k)\|_\infty \leq \epsilon_P,
\end{equation}
where \(\epsilon_D > 0\) and \(\epsilon_P > 0\) are user-defined tolerances.
Using the infinity norm in~\eqref{eq:overall-stopping-criterion} is common for nonlinear optimization solvers \citep{gill-murray-saunders-2002, wachter-biegler-2006}.
To monitor the accuracy of the approximate solution \((x_k, y_k)\) of~\eqref{eq:penalty-nlp} we use a sequence of primal and dual tolerances \(\epsilon_P^k > 0\) and \(\epsilon_D^k > 0\).
In view of~\eqref{eq:outer-kkt}, we impose that the pair \((x_k, y_k)\) satisfy
\begin{equation}%
  \label{eq:accuracy-condition}
  \| \nabla f(x_k) + { J(x_k) }^T y_k \|_\infty \leq \epsilon_D^k.
\end{equation}

If \(\|c(x_k)\|_\infty > \epsilon_P^k\), we increase \(\tau_k\) according to
\begin{equation}%
  \label{eq:penalty-update}
  \tau_{k + 1} = \max(\tau_k + \tau_{\min}, \|y_k\|_1),
\end{equation}
where \(\tau_{\min} > 0\) is a minimum increase for the penalty parameter.
We choose the \(1\)-norm in~\eqref{eq:penalty-update} to overestimate \(\|\bar y\|_2\), using the fact that \(\|\cdot\|_2 \leq \|\cdot\|_1\).

If instead \(\|c(x_k)\|_\infty \leq \epsilon_P^k\), then, since \((x_k, y_k)\) also satisfy~\eqref{eq:accuracy-condition}, \((x_k, y_k)\) satisfy~\eqref{eq:overall-stopping-criterion} for the tolerances \((\epsilon_D^k, \epsilon_P^k)\).
Therefore, we tighten the accuracy on \((x_k, y_k)\) by updating
\begin{equation}%
  \label{eq:accuracy-update}
  \epsilon_D^{k + 1} = \max(\mu_{\epsilon}\epsilon_D^k, \epsilon_D) \quad \text{and} \quad \epsilon_P^{k + 1} = \max(\mu_{\epsilon}\epsilon_P^k, \epsilon_P),
\end{equation}
where \(\mu_{\epsilon} \in (0, 1)\).

We reset the accuracy for the subproblem, \(\epsilon_D^{k + 1} = \epsilon_D^1\), each time we increase \(\tau_k\), because \(x_k\) may not be an accurate estimate of the solution of the new penalized problem~\eqref{eq:penalty-nlp}.
In addition, we improve the initial guess for the next penalized problem~\eqref{eq:penalty-nlp} by applying an extrapolation technique to \((x_k, y_k)\).
Our approach is similar to extrapolation techniques used in interior-point methods \citep{gould-orban-sartenaer-toint-2001}.

Beyond approximately solving~\eqref{eq:penalty-nlp}, the ``outer loop'' iterations also detect the local infeasibility of~\eqref{eq:nlp} at \(x_k\) by computing an optimality measure for
\begin{equation}%
  \label{eq:infeasibility-problem}
  \minimize{x \in \R^n} \ \|c(x)\|_2.
\end{equation}
When we detect a local solution of~\eqref{eq:infeasibility-problem} with a strictly positive objective value, we declare~\eqref{eq:nlp} to be \emph{locally} infeasible.

\Cref{alg:exactpen-outer-loop} summarizes the ``outer loop'' of our exact penalty method.

\begin{algorithm}[ht]%
  \caption{%
    \label{alg:exactpen-outer-loop}
    Outer Loop.
  }
  \begin{algorithmic}[1]%

    \Require{}%
    Starting points \((x_0, y_0)\);
    initial penalty parameter \(\tau_{\mathrm{init}} > 0\);
    constants \(\epsilon_D^1 \geq \epsilon_D > 0\);
    \(\epsilon_P^1 \geq \epsilon_P > 0\);
    \(\tau_{\min} > 0\);
    \(\mu_{\epsilon} \in (0, 1)\).

    \State{}%
    Check~\eqref{eq:overall-stopping-criterion} for \((x_0, y_0)\).
    If it is satisfied, \textbf{stop} and declare \textbf{first-order stationarity}.

    \State{}%
    \label{line:init-tau}
    Set \(\tau_1 \coloneqq \max(\tau_{\mathrm{init}}, \|y_0\|_1)\).

    \For{\(k = 1, \ldots\)}%

    \State{}%
    Compute \((x_{k}, y_{k})\), an approximate solution of~\eqref{eq:penalty-nlp} that satisfies~\eqref{eq:accuracy-condition}, started from a well-chosen initial guess.

    \If{\(\|c(x_k)\|_{\infty} > \epsilon_P^k\)}%
    \Comment{Increase penalty parameter.}

    \State{}%
    \label{line:penalty-update}%
    Update \(\tau_{k+1}\) according to~\eqref{eq:penalty-update}.

    \State{}%
    Update \((x_k, y_k)\) with the extrapolation technique.

    \State%
    Reset \(\epsilon_D^{k + 1} = \epsilon_D^1\) and set \(\epsilon_P^{k + 1} = \epsilon_P^{k}\).

    \Else%
    \Comment{Tighten tolerance.}
    \State%
    Set \(\tau_{k+1} \coloneqq \tau_k\).

    \State%
    Update \(\epsilon_D^{k + 1}\) and \(\epsilon_P^{k + 1}\) with~\eqref{eq:accuracy-update}.

    \EndIf%

    \State{}%
    Check~\eqref{eq:overall-stopping-criterion} for \((x_k, y_k)\).
    If satisfied, \textbf{stop} and declare \textbf{first-order stationarity}.

    \State{}%
    Estimate local infeasibility.

    \EndFor%
  \end{algorithmic}
\end{algorithm}

\subsection{Extrapolation Technique}%
\label{sec:extrapolation-technique}

Let \((x_k, y_k)\) be an approximate solution of~\eqref{eq:penalty-nlp} for a given penalty parameter \(\tau_k\).
In \Cref{alg:exactpen-outer-loop}, when increasing the penalty parameter to \(\tau_{k+1}\), we use an extrapolation technique to improve the initial guess for the next penalized problem~\eqref{eq:penalty-nlp}.
Indeed, \((x_k, y_k)\) may not be an accurate estimate of the solution of~\eqref{eq:penalty-nlp} for \(\tau \coloneqq \tau_{k+1}\).
To do so, we apply the implicit function theorem to the first-order optimality conditions of~\eqref{eq:penalty-nlp} derived in \Cref{prop:first-order-penalized}, whose equality conditions are smooth even though~\eqref{eq:penalty-nlp} is nonsmooth.

We collect the equality conditions in~\eqref{eq:outer-kkt} and~\eqref{eq:inner-kkt} into the mapping
\[
  F(x, y, \alpha, \tau) \coloneqq
  \begin{pmatrix}
    \nabla f(x) + { J(x) }^T y \\
    c(x) - \alpha y            \\
    \alpha(\|y\|_2^2 - \tau^2)
  \end{pmatrix}.
\]
Thus, if \(z^* \coloneqq (x^*, y^*, \alpha^*)\) satisfies the first-order optimality conditions of~\eqref{eq:penalty-nlp} at \(\tau^*\), then
\(
F(z^*, \tau^*) = 0.
\)
Assume that the Jacobian of \(F\) with respect to \((x, y, \alpha)\) is nonsingular at \((z^*, \tau^*)\).
Then, the implicit function theorem implies existence of a \(\mathcal{C}^1\) path \(\tau \to z_\tau \coloneqq (x_\tau, y_\tau, \alpha_\tau)\), defined in a neighborhood of \(\tau^*\), such that
\(
z_{\tau^*} = z^*,
\)
and such that the derivative of the path \(\dot z_{\tau^*}\) satisfies
\begin{equation}%
  \label{eq:ift-system}
  D_{z} F(z^*, \tau^*)
  \dot z_{\tau^*}
  =
  - D_\tau F(z^*, \tau^*).
\end{equation}
For completeness, we expand the Jacobians in~\eqref{eq:ift-system} as
\[
  D_z F(z, \tau)
  =
  \begin{pmatrix}
    \nabla_{xx}^2 \mathcal{L}(x, y) & {J(x)}^T     & 0                  \\
    J(x)                            & -\alpha I    & -y                 \\
    0                               & 2 \alpha y^T & \|y\|_2^2 - \tau^2
  \end{pmatrix},
  \quad
  \textup{and}
  \quad
  D_\tau F(z, \tau)
  =
  \begin{pmatrix}
    0 \\
    0 \\
    -2 \alpha \tau
  \end{pmatrix}.
\]
At \Cref{line:penalty-update} of \Cref{alg:exactpen-outer-loop}, given a new penalty parameter \(\tau_{k+1}\), we compute
\[
  x_{k, \mathrm{extr}}
  \coloneqq
  x_k
  +
  (\tau_{k+1} - \tau_k)\dot x_k,
  \quad
  \text{and}
  \quad
  y_{k, \mathrm{extr}}
  \coloneqq
  y_k
  +
  (\tau_{k+1} - \tau_k)\dot y_k,
\]
where \(\dot z_k = (\dot x_k, \dot y_k, \dot \alpha_k)\) is the solution of~\eqref{eq:ift-system} for \(\tau^* = \tau_k\) and \(z^* = (x_k, y_k, \alpha_k)\), and \(\alpha_k\) is defined at the end of \Cref{sec:step-computation}.
We start the ``inner loop,'' described in \Cref{sec:penalized-subproblem}, with \((x_{k, \mathrm{extr}}, y_{k, \mathrm{extr}})\).

\subsection{Infeasibility Detection}%
\label{sec:infeasibility-detection}

To detect local infeasibility of~\eqref{eq:nlp}, we use an optimality measure for~\eqref{eq:infeasibility-problem} adapted from \citet{diouane-gollier-orban-2026}.
Specifically, we use
\begin{equation}
  \label{eq:infeasibility-optimality-measure}
  \sqrt{\nu_k^{-1} \theta(x_k) } \coloneqq \sqrt{\nu_k^{-1}}\sqrt{\|c(x_k)\|_2 - \|c(x_k) + {J(x_k)} s_{k, \mathrm{feas}}\|_2},
\end{equation}
where
\begin{equation}%
  \label{eq:prox-feasibility}
  s_{k, \mathrm{feas}} \in \argmin{s \in \R^n} \|c(x_k) + J(x_k) s\|_2 + \tfrac{1}{2 \nu_k} \|s\|_2^2,
\end{equation}
and \(\nu_k > 0\) is a regularization parameter.
\begin{proposition}%
  \label{prop:infeasibility-optimality-measure}
  For any \(\nu_k > 0\), if \(\sqrt{\nu_k^{-1} \theta(x_k)} = 0\), then \(x_k\) is first-order stationary for~\eqref{eq:infeasibility-problem}.
\end{proposition}
\begin{proof}%
  This is a consequence of the discussion following equation~(4.8) in \citep[Section 4]{diouane-gollier-orban-2026}.
\end{proof}
To compute a solution of~\eqref{eq:prox-feasibility}, we use the following result.

\begin{proposition}%
  \label{prop:infeasibility-computation}
  Assume that \(c(x_k) \neq 0\).
  The unique solution \(s_{k, \mathrm{feas}}\) to~\eqref{eq:prox-feasibility} satisfies
  \begin{equation}%
    \label{eq:prox-feasibility-system}
    \begin{pmatrix}
      I      & {J(x_k)}^T \\
      J(x_k) & - I
    \end{pmatrix}
    \begin{pmatrix}
      s_{k, \mathrm{feas}} \\
      y_{k, \mathrm{feas}}
    \end{pmatrix}
    =
    -
    \begin{pmatrix}
      0 \\
      c(x_k)
    \end{pmatrix}.
  \end{equation}
  The regularization parameter \(\nu_k > 0\) of~\eqref{eq:prox-feasibility} corresponding to \(s_{k, \mathrm{feas}}\) is \(\nu_k \coloneqq \|y_{k, \mathrm{feas}}\|_2\).
\end{proposition}
\begin{proof}%
  Using the same proof technique as in \Cref{prop:first-order-penalized} with \(\tau = 1\), and since~\eqref{eq:prox-feasibility} is strongly convex, its unique solution is characterized by the necessary and sufficient optimality conditions
  \[
    \nu_k^{-1} s + { J(x_k) }^T y = 0, \quad
    c(x_k) + J(x_k) s - \alpha y = 0, \quad
    \alpha (\|y\|_2^2 - 1) = 0, \quad
    \|y\|_2 \leq 1, \quad
    \alpha \geq 0,
  \]
  which are satisfied for \((s, y, \alpha) = (s_{k, \mathrm{feas}}, \nu_k^{-1}y_{k, \mathrm{feas}}, \nu_k)\), with \(\nu_k = \|y_{k, \mathrm{feas}}\|_2\).
  To see this, the first row of~\eqref{eq:prox-feasibility-system} yields
  \[
    \nu_k^{-1} s + { J(x_k) }^T y = \nu_k^{-1} s_{k, \mathrm{feas}} + \nu_k^{-1} {J(x_k)}^T y_{k, \mathrm{feas}} = 0.
  \]
  Moreover, the second row of~\eqref{eq:prox-feasibility-system} gives
  \[
    c(x_k) + J(x_k) s_{k, \mathrm{feas}} - \nu_k \left( \nu_k^{-1} y_{k, \mathrm{feas}} \right)
    = c(x_k) + J(x_k) s_{k, \mathrm{feas}} - y_{k, \mathrm{feas}}
    = 0.
  \]
  Finally, since \(\nu_k = \|y_{k, \mathrm{feas}}\|_2\),
  \[
    \|\nu_k^{-1} y_{k,\mathrm{feas}}\|_2 = \nu_k^{-1} \|y_{k,\mathrm{feas}}\|_2 = \nu_k^{-1} \nu_k = 1 \leq 1,
  \]
  so that \(\alpha (\|y\|_2^2 - 1) = \nu_k (1 - 1) = 0\) and \(\alpha = \nu_k > 0\), which completes the verification that \((s,y,\alpha) = (s_{k,\mathrm{feas}}, \nu_k^{-1} y_{k,\mathrm{feas}}, \nu_k)\) satisfies the optimality conditions of~\eqref{eq:prox-feasibility}.

  Note that \(\nu_k = \|y_{k, \mathrm{feas}}\|_2 > 0\) because \(y_{k, \mathrm{feas}} = 0\) implies from the first row of~\eqref{eq:prox-feasibility-system} that \(s_{k, \mathrm{feas}} = 0\), which then gives \(y_{k, \mathrm{feas}} = c(x_k) \neq 0\), a contradiction.
\end{proof}
In contrast to \citet{diouane-gollier-orban-2026}, where \(\nu_k = 1\) is fixed before solving~\eqref{eq:prox-feasibility} with \citep[Corollary~5.3]{diouane-gollier-orban-2026}, \(\nu_k\) is here determined by the solution of~\eqref{eq:prox-feasibility-system}.

Following \Cref{prop:infeasibility-computation}, we compute~\eqref{eq:infeasibility-optimality-measure} by solving~\eqref{eq:prox-feasibility-system} for \((s_{k, \mathrm{feas}}, y_{k, \mathrm{feas}})\) and setting \(\nu_k = \|y_{k, \mathrm{feas}}\|_2\).
We wish to declare local infeasibility when \(\sqrt{\nu_{k}^{-1}\theta(x_k)}\) is small.
Specifically, we impose that \(\sqrt{\nu_k^{-1}\theta(x_k)} < \epsilon_I \|c(x_k)\|_2\), where \(\epsilon_I > 0\) is an algorithmic constant, to declare local infeasibility.
To guarantee that \(x_k\) is close to a stationary point of~\eqref{eq:infeasibility-problem}, we additionally require \(\sqrt{\nu_k^{-1}\theta(x_k)}\) to be smaller than the current primal tolerance \(\epsilon_P^k\).
In summary, we declare the problem~\eqref{eq:nlp} to be locally infeasible when
\[
  \sqrt{\nu_k^{-1}\theta(x_k)} < \epsilon_I \|c(x_k)\|_2,\quad
  \|c(x_k)\|_2 > \epsilon_P
  \quad \textup{and} \quad
  \sqrt{\nu_k^{-1}\theta(x_k)} < \epsilon_P^k,
\]
where \(\epsilon_P\), \(\epsilon_P^k\) are the tolerances described in \Cref{sec:exact-penalty-method}.

\subsection{Solution of the Penalized Problem}%
\label{sec:penalized-subproblem}

To solve~\eqref{eq:penalty-nlp} for fixed \(\tau_k\), we use an iterative algorithm inspired by \citep[Algorithm 4.1]{diouane-habiboullah-orban-2026}.
We refer to this process as the ``inner loop'' of our algorithm.
Let \(l\) be the iteration counter of the ``inner loop.''
The method computes a step \(s_l\) by approximately solving a regularized quadratic model of~\eqref{eq:penalty-nlp} around the current iterate \(x_l\) which is,
\begin{equation}%
  \label{eq:quadratic-subproblem}
  \minimize{s \in \R^n} \ q(s; x_l, \sigma_l) \coloneqq f_l + {\nabla f_l}^T s + \tfrac{1}{2} s^T (H_l + \sigma_l I) s + \tau_k \|c_l + J_l s\|_2,
\end{equation}
where \(\sigma_l > 0\) is a quadratic regularization parameter and \(H_l\) is a symmetric matrix.

We accept a step that satisfies
\begin{equation}%
  \label{eq:armijo}
  \rho_l \coloneqq
  \frac{
  f_l + \tau_k \|c_l\|_2 - \left( f(x_l + s_l) + \tau_k \|c(x_l + s_l)\|_2 \right)
  }{
  f_l + \tau_k \|c_l\|_2 - (f_l + {\nabla f_l}^T s_l + \tau_k\|c_l + J_l s_l\|_2)
  } \geq \eta_1,
\end{equation}
for some constant \(\eta_1 \in (0, 1)\), in which case we set \(x_{l+1} \coloneqq x_l + s_l\).
If \(\rho_l < \eta_1\), then \(x_{l+1} \coloneqq x_l\) and the step is rejected.
When we reject a step, we increase the regularization parameter \(\sigma_{l+1} > \sigma_l\).
When \(\rho_l \geq \eta_2 \in [\eta_1, 1)\), we decrease \(\sigma_{l+1} \leq \sigma_l\).

As we shall see in \Cref{sec:convergence-analysis}, under a Lipschitz assumption on \(\nabla f\) and \(J\), steps are guaranteed to be accepted when \(\sigma_l\) is sufficiently large.
This fact ensures global convergence of the method.
On the other hand, though not proven, we observe fast local convergence of our method close to solutions of~\eqref{eq:penalty-nlp} when \(\sigma_l \approx 0\).
For this reason, the parameter \(\sigma_l\) plays a crucial role in the performance of the method in practice.
The basic update rule that guarantees the convergence is given in \Cref{alg:exactpen-inner-loop} and accelerating heuristics are provided in \Cref{sec:regularization-updates}.

We assume that the computation of the step \(s_l\) yields an estimate \(y_{l+1}\) of the Lagrange multipliers, as formalized below.
\begin{stepassumption}%
  \label{stepassumption:step-s-stationarity}
  For each \(l \in \N\), the step \(s_l\) satisfies
  \begin{equation}%
    \label{eq:step-s-stationarity}
    (H_l + \sigma_l I) s_l = - \nabla f_l - { J_l }^T y_{l+1},
  \end{equation}
  for some Lagrange multiplier estimate \(y_{l+1} \in \R^m\).
\end{stepassumption}

For the method to converge, we make the following additional assumption on the step \(s_l\),
\begin{stepassumption}%
  \label{stepassumption:decrease}
  For each \(l \in \N\), the step \(s_l\) satisfies the simple decrease condition
  \begin{equation}%
    \label{eq:simple-decrease}
    q(s_l; x_l, \sigma_l) \leq q(0; x_l, \sigma_l),
  \end{equation}
  and, for some constant \(\eta_C > 0\) independent of \(l\), the Cauchy decrease condition
  \begin{equation}%
    \label{eq:cauchy-decrease}
    f_l + \tau_k \|c_l\|_2 - f_l - {\nabla f_l}^T s_l - \tau_k \|c_l + J_l s_l\|_2 \geq \tfrac{1}{2} \eta_C \min(1, \|\nabla f_l + J_l^T y_{l+1}\|_{\infty}^2).
  \end{equation}
\end{stepassumption}

As a consequence of \Cref{stepassumption:step-s-stationarity}, our solver produces iterates \((x_l, y_l)\), which we use to define \(H_l\) in~\eqref{eq:quadratic-subproblem} as either the Hessian of the Lagrangian \(\nabla^2 \mathcal{L}(x_l, y_l)\) or a quasi-Newton approximation thereof.
When we present our step computation procedure in \Cref{sec:step-computation}, we will show that our implementation satisfies \Cref{stepassumption:step-s-stationarity,stepassumption:decrease}.
\Cref{stepassumption:decrease} stands as a significant improvement over the penalty method in \citep{diouane-gollier-orban-2026}, which assumed that the step \(s_l\) is an exact solution of~\eqref{eq:quadratic-subproblem}.

\Cref{alg:exactpen-inner-loop} summarizes the ``inner loop'' of our exact penalty method.

\begin{algorithm}[ht]%
  \caption{%
    \label{alg:exactpen-inner-loop}
    Inner Loop.
  }
  \begin{algorithmic}[1]%

    \Require{}%
    Starting points \((x_0, y_0)\);
    initial regularization parameter \(\sigma_0 > 0\);
    penalty parameter \(\tau_k > 0\);
    constants \(\epsilon_D^k > 0\);
    \(\eta_1 \in (0, 1)\);
    \(\eta_2 \in [\eta_1, 1)\);
      \(\gamma_1 > 1\);
      \(\gamma_2 \geq \gamma_1\);
      \(\gamma_3 \in (0, 1]\);

    \For{\(l = 0, 1, \ldots\)}%

    \State%
    Choose a symmetric matrix \(H_l \approx \nabla^2 \mathcal{L}(x_l, y_l)\).

    \State%
    Compute a step \(s_l\) and an estimate of the Lagrange multipliers \(y_{l+1}\) that satisfy~\Cref{stepassumption:step-s-stationarity,stepassumption:decrease}.

    \State%
    Check~\eqref{eq:accuracy-condition} for \((x_l, y_{l+1})\). If satisfied, \textbf{stop} and return \((x_l, y_{l+1})\).

    \State%
    Compute the ratio
    \begin{equation*}%
      \rho_l \coloneqq \frac{f_l + \tau_k \|c_l\|_2 - \left( f(x_l + s_l) + \tau_k \|c(x_l + s_l)\|_2 \right) }{f_l + \tau_k \|c_l\|_2 - (f_l + {\nabla f_l}^T s_l + \tau_k\|c_l + J_l s_l\|_2)}.
    \end{equation*}

    \State%
    If \(\rho_l \geq \eta_1\), set \(x_{l+1} \coloneqq x_l + s_l\).
    Otherwise, set \(x_{l+1} \coloneqq x_l\).

    \State%
    \label{algline:regularization-update}
    Update the regularization parameter according to
    \begin{equation*}
      \sigma_{l+1} \in
      \begin{cases}
        \left[ \gamma_3\sigma_l, \sigma_l \right]         & \text{if} \ \rho_l \geq \eta_2,          \\
        \sigma_l                                          & \text{if} \ \eta_1 \leq \rho_l < \eta_2, \\
        \left[ \gamma_1\sigma_l, \gamma_2\sigma_l \right] & \text{if} \ \rho_l < \eta_1.
      \end{cases}
    \end{equation*}

    \EndFor%
  \end{algorithmic}
\end{algorithm}

\section{Convergence Analysis}%
\label{sec:convergence-analysis}

The aim of this section is to show that \Cref{alg:exactpen-outer-loop} produces iterates \((x_k, y_k)\) that eventually satisfy~\eqref{eq:overall-stopping-criterion}.
Our approach can be summarized as follows.
First, we show that \Cref{alg:exactpen-inner-loop} brings
\begin{equation}%
  \label{eq:optimality-measure}
  \| \nabla f(x_l) + {J(x_l)}^T y_{l+1} \|_{\infty},
\end{equation}
below a given tolerance \(\epsilon > 0\) in \(\mathcal{O}(\epsilon^{-2})\) iterations.
Then, we show convergence of \Cref{alg:exactpen-outer-loop}.

\subsection{Convergence of the Inner Loop}%
\label{sec:inner-loop-convergence}

We now prove a complexity result for \Cref{alg:exactpen-inner-loop} with respect to~\eqref{eq:optimality-measure}.
Our working assumptions are the following.
\begin{modelassumption}%
  \label{modelassumption:lipschitz}
  There exists \(\kappa_m > 0\) such that for all \(l \in \N\) and \(k \in \N\),
  \begin{equation*}%
    | f(x_l + s_l) + \tau_k \|c(x_l + s_l)\|_2 - (f_l + {\nabla f_l}^T s_l + \tau_k \|c_l + J_l s_l\|_2 )| \leq \kappa_m \|s_l\|_2^2.
  \end{equation*}
\end{modelassumption}

\begin{modelassumption}%
  \label{modelassumption:bounded-hessian}
  There exists \(\kappa_H > 0\) such that for all \(l \in \N\) and \(k \in \N\),
  \begin{equation*}%
    \|H_l\|_2 \leq \kappa_H.
  \end{equation*}
\end{modelassumption}

Following the discussion in \citep[Section 4]{diouane-gollier-orban-2026}, \Cref{modelassumption:lipschitz} is satisfied when \(\nabla f\) and \(J\) are Lipschitz continuous and the penalty parameter \(\tau_k\) is bounded.

The following lemma is critical in showing global convergence of the method.

\begin{lemma}%
  \label{lemma:sigma_succ}
  Let condition~\eqref{eq:simple-decrease} of \Cref{stepassumption:decrease} and \Cref{modelassumption:lipschitz,modelassumption:bounded-hessian} be satisfied.
  Define
  \begin{equation}%
    \label{eq:sigma_succ}
    \sigma_{\mathrm{succ}} \coloneqq \kappa_H + \tfrac{2 \kappa_m }{ (1 - \eta_2) }.
  \end{equation}
  If, at iteration \(l\) of \Cref{alg:exactpen-inner-loop}, \(s_l \neq 0\) and \(\sigma_l \geq \sigma_{\mathrm{succ}}\), then \(\rho_l \geq \eta_2\).
\end{lemma}
\begin{proof}%
  From~\eqref{eq:sigma_succ}, \(\sigma_l \geq \sigma_{\mathrm{succ}} > \kappa_H\), so that \Cref{modelassumption:bounded-hessian} implies that \(H_l + \sigma_l I\) is positive definite.
  Additionally, since \(s_l\) satisfies~\eqref{eq:simple-decrease},
  \begin{equation}%
    \label{eq:cp-model-decrease}
    \tau_k \|c_l\|_2 - {\nabla f_l}^T s_l - \tau_k \|c_l + J_l s_l\|_2 \geq \tfrac{1}{2} s_l^T (H_l + \sigma_l I) s_l.
  \end{equation}
  Then, the definition of \(\rho_l\) in \Cref{alg:exactpen-inner-loop},~\eqref{eq:cp-model-decrease}, and \Cref{modelassumption:lipschitz} imply that
  \begin{align*}%
    \label{eq:rho-j-prop}
    |\rho_l - 1| & =
    \left|
    \frac{
    f(x_l + s_l) + \tau_k \|c(x_l + s_l)\|_2 - f_l - {\nabla f_l}^T s_l - \tau_k \|c_l + J_l s_l\|_2
    }{
    \tau_k \|c_l\|_2 - {\nabla f_l}^T s_l - \tau_k \|c_l + J_l s_l\|_2
    }
    \right|             \\
                 & \leq
    \frac{
      2 \kappa_m \|s_l\|_2^2
    }{
      {s_l}^T (H_l + \sigma_l I) s_l
    }                   \\
                 & \leq
    \frac{
      2 \kappa_m
    }{
      \lambda_{\min}(H_l + \sigma_l I)
    }.
  \end{align*}
  Now, from~\eqref{eq:sigma_succ},
  \begin{equation*}%
    \lambda_{\min}(H_l + \sigma_l I) = \lambda_{\min}(H_l) + \sigma_l \geq \sigma_l - \kappa_H \geq \frac{2 \kappa_m}{1 - \eta_2},
  \end{equation*}
  which gives
  \(
  |\rho_l - 1| \leq 1 - \eta_2.
  \)
\end{proof}

The mechanism of \Cref{alg:exactpen-inner-loop} and \Cref{lemma:sigma_succ} imply that
\begin{equation}%
  \label{eq:sigma-max}
  \sigma_l \leq \sigma_{\max} \coloneqq \max\{\sigma_0, \gamma_2\sigma_{\mathrm{succ}}\}, \quad \forall l \in \N.
\end{equation}
As a consequence of~\eqref{eq:sigma-max}, we establish the following connection between the two conditions of \Cref{stepassumption:decrease}.
\begin{lemma}%
  \label{lemma:simple-cauchy-decrease}
  Let \Cref{stepassumption:step-s-stationarity} and condition~\eqref{eq:simple-decrease} of \Cref{stepassumption:decrease} be satisfied and \Cref{modelassumption:lipschitz,modelassumption:bounded-hessian} be satisfied.
  Then, for each \(l \in \N\), if \(H_l + \sigma_l I\) is positive definite, \(s_l\) satisfies~\eqref{eq:cauchy-decrease} with \(\eta_C = {(\kappa_H + \sigma_{\max})}^{-1}\).
\end{lemma}
\begin{proof}%
  From~\eqref{eq:simple-decrease}, \(s_l\) satisfies~\eqref{eq:cp-model-decrease}.
  Moreover, under~\eqref{eq:simple-decrease} and \Cref{modelassumption:lipschitz,modelassumption:bounded-hessian}, \Cref{lemma:sigma_succ} implies~\eqref{eq:sigma-max}.
  Then, from \(H_l + \sigma_l I \succ 0\), \Cref{modelassumption:bounded-hessian} and~\eqref{eq:sigma-max},
  \begin{align*}
    \tau_k \|c_l\|_2 - {\nabla f_l}^T s_l - \tau_k \|c_l + J_l s_l\|_2
     & \geq
    \tfrac{1}{2} s_l^T (H_l + \sigma_l I) s_l                                              \\
     & =
    \tfrac{1}{2} s_l^T (H_l + \sigma_l I) {(H_l + \sigma_l I)}^{-1} (H_l + \sigma_l I) s_l \\
     & \geq
    \tfrac{1}{2} \lambda_{\min}({(H_l + \sigma_l I)}^{-1}) \|(H_l + \sigma_l I) s_l\|_2^2  \\
     & \geq
    \tfrac{1}{2} {(\kappa_H + \sigma_{\max})}^{-1} \|(H_l + \sigma_l I) s_l\|_2^2.
  \end{align*}
  The conclusion follows by injecting~\eqref{eq:step-s-stationarity} into the last inequality and using the fact that \(\|\cdot\|_2 \geq \|\cdot\|_\infty\).
\end{proof}

Moving on to the complexity analysis, given a user-defined tolerance \(\epsilon \leq 1\), we seek an upper bound on \(l_{\epsilon} \coloneqq \min\{ l \in \N \mid \| \nabla f_l + J_l^T y_{l+1} \|_{\infty} \leq \epsilon \} = |\mathcal{S}_k(\epsilon)| + |\mathcal{U}_k(\epsilon)| + 1\), where
\begin{equation*}%
  \mathcal{S}_k(\epsilon) \coloneqq \{ l \in \N \mid \rho_l \geq \eta_1 \ \text{and} \ l < l_{\epsilon}\}, \quad
  \mathcal{U}_k(\epsilon) \coloneqq \{ l \in \N \mid \rho_l < \eta_1 \ \text{and} \ l < l_{\epsilon}\}.
\end{equation*}
We use subscripts \(k\) to emphasize that the complexity analysis is for the inner loop of outer loop iteration \(k\) of \Cref{alg:exactpen-outer-loop}.

\begin{lemma}%
  \label{lemma:S-complexity}
  Let condition~\eqref{eq:cauchy-decrease} of \Cref{stepassumption:decrease}, \Cref{stepassumption:step-s-stationarity}, \Cref{modelassumption:lipschitz,modelassumption:bounded-hessian} be satisfied.
  Assume that \(x \to f(x) + \tau_k \|c(x)\|_2\) is bounded below by \({(f + \tau_k \|c\|_2)}_{\mathrm{low}}\).
  Then,
  \begin{equation}%
    \label{eq:S-bound}
    |\mathcal{S}_k(\epsilon)| \leq
    \frac{
    f_0 + \tau_k \|c_0\|_2 - {(f + \tau_k \|c\|_2)}_{\mathrm{low}}
    }{
    \tfrac{1}{2} \eta_1 \eta_C
    }
    \epsilon^{-2}.
  \end{equation}
\end{lemma}
\begin{proof}%
  Let \(l \in \mathcal{S}_k(\epsilon)\).
  From~\eqref{eq:armijo} and~\eqref{eq:cauchy-decrease},
  \begin{align*}
    f_l + \tau_k \|c_l\|_2 - f(x_l + s_l) - \tau_k \|c(x_l + s_l) \|_2
     & \geq
    \eta_1 (\tau_k \|c_l\|_2 - {\nabla f_l}^T s_l - \tau_k \|c_l + J_l s_l\|_2) \\
     & \geq
    \tfrac{1}{2} \eta_1 \eta_C   \min(1, \|{\nabla f_l} + J_l^T y_{l+1}\|_\infty^2).
  \end{align*}
  Because \(1 \geq \epsilon^2\) by definition and \(\|\nabla f_l + J_l^T y_{l+1}\|_\infty^2 > \epsilon^2\) following that \(l < l_{\epsilon}\), the above simplifies to
  \begin{align*}
    f_l + \tau_k \|c_l\|_2 - f(x_l + s_l) - \tau_k \|c(x_l + s_l) \|_2
     & \geq
    \tfrac{1}{2} \eta_1 \eta_C \epsilon^2.
  \end{align*}
  The remainder of the proof is classic and identical to, e.g., \citep[Lemma 3.6]{aravkin-baraldi-orban-2022}.
\end{proof}
\begin{lemma}%
  \label{lemma:U-complexity}
  Under the assumptions of \Cref{lemma:sigma_succ,lemma:S-complexity},
  \begin{equation}%
    \label{eq:U-bound}
    |\mathcal{U}_k(\epsilon)| \leq \log_{\gamma_1}(\sigma_{\max}/\sigma_0) + | \mathcal{S}_k(\epsilon) || \log_{\gamma_1}(\gamma_3) |.
  \end{equation}
\end{lemma}
\begin{proof}
  The proof is similar to \citep[Lemma 3.7]{aravkin-baraldi-orban-2022} and is omitted.
\end{proof}

\begin{theorem}%
  \label{thm:inner-loop-complexity}
  Let \Cref{stepassumption:decrease}, \Cref{stepassumption:step-s-stationarity}, and \Cref{modelassumption:bounded-hessian,modelassumption:lipschitz} hold and assume that \(x \to f(x) + \tau_k \|c(x)\|_2\) is bounded below by \({(f + \tau_k \|c\|_2)}_{\mathrm{low}}\).
  Then,
  \begin{equation}%
    \label{eq:inner-loop-bound}
    l_{\epsilon} \leq
    \left( 1 + |\log_{\gamma_1}(\gamma_3)| \right)
    \frac{
    f_0 + \tau_k \|c_0\|_2 - {(f + \tau_k \|c\|_2)}_{\mathrm{low}}
    }{
    \tfrac{1}{2} \eta_1 \eta_C
    }
    \epsilon^{-2}
    + \log_{\gamma_1}(\sigma_{\max}/\sigma_0) + 1.
  \end{equation}
  In particular, \(l_{\epsilon} = \mathcal{O}(\epsilon^{-2})\).
\end{theorem}
\begin{proof}%
  By definition, \(l_{\epsilon} = |\mathcal{S}_k(\epsilon)| + |\mathcal{U}_k(\epsilon)| + 1\).
  The result follows by bounding \(|\mathcal{U}_k(\epsilon)|\) with~\eqref{eq:U-bound} and then \(|\mathcal{S}_k(\epsilon)|\) with~\eqref{eq:S-bound}.
\end{proof}

\subsection{Convergence of the Outer Loop}%
\label{sec:outer-loop-convergence}

\Cref{thm:inner-loop-complexity} allows us to establish that the iterates \((x_k, y_k)\) generated by \Cref{alg:exactpen-outer-loop} are well defined and satisfy~\eqref{eq:accuracy-condition}.

\begin{corollary}%
  \label{cor:outer-iterates-well-defined}
  Let \Cref{stepassumption:decrease,stepassumption:step-s-stationarity} and \Cref{modelassumption:lipschitz,modelassumption:bounded-hessian} hold and assume that \(x \to f(x) + \tau_1 \|c(x)\|_2\) is bounded below by \({(f + \tau_1 \|c\|_2)}_{\mathrm{low}}\).
  Then, for each \(k \in \N\), the iterates \((x_k, y_k)\) of \Cref{alg:exactpen-outer-loop} generated by \Cref{alg:exactpen-inner-loop} are well defined and satisfy~\eqref{eq:accuracy-condition}.
\end{corollary}
\begin{proof}%
  If \(x \to f(x) + \tau_1 \|c(x)\|_2\) is bounded below, then \(x \to f(x) + \tau_k \|c(x)\|_2\) is also bounded below for all \(k \in \N\) since \(\tau_k \geq \tau_1\).
  Therefore, because \Cref{stepassumption:decrease,stepassumption:step-s-stationarity} and \Cref{modelassumption:lipschitz,modelassumption:bounded-hessian} hold, \Cref{thm:inner-loop-complexity} applies for each \(k \in \N\), and \Cref{alg:exactpen-inner-loop} brings~\eqref{eq:optimality-measure} below \(\epsilon_D^k\) in \(\mathcal{O}({(\epsilon_D^k)}^{-2})\) iterations.
  Then, since \Cref{stepassumption:step-s-stationarity} holds, the pair \((x_k, y_k)\) satisfies~\eqref{eq:accuracy-condition}.
\end{proof}

The following lemma proves the convergence of \Cref{alg:exactpen-outer-loop} to first-order stationary points of~\eqref{eq:nlp}.
The approach here is similar to \citep[Theorem 4.5]{diouane-gollier-orban-2026} with the exception that we assume that the penalty parameter is bounded.
Proving the convergence in the unbounded case is beyond the scope of this paper.
For simplicity, we assume that \(\epsilon_P = \epsilon_D = \epsilon\).
We seek an upper bound on \(k_\epsilon \coloneqq \min\{ k \in \N \mid \|\nabla f_k + J_k^T y_k\|_\infty \leq \epsilon \ \text{and} \ \|c_k\|_\infty \leq \epsilon \}\).
\begin{lemma}%
  \label{lemma:outer-loop-conv}
  Let the assumptions of \Cref{cor:outer-iterates-well-defined} hold and assume that there is a \(\bar \tau\), independent of \(\epsilon\) such that \(\|c(x_k)\|_\infty \leq \epsilon\) whenever \(\tau_k \geq \bar \tau\).
  Then,
  \[
    k_\epsilon \leq
    \left\lceil \frac{\bar \tau - \tau_{\mathrm{init}}}{\tau_{\min}} \right\rceil +
    \left\lceil \log_{\mu_\epsilon}(\epsilon / \epsilon_P^1) \right\rceil +
    \left( \left\lceil \frac{\bar \tau - \tau_{\mathrm{init}}}{\tau_{\min}} \right\rceil + 1 \right)
    \left\lceil \log_{\mu_\epsilon}(\epsilon / \epsilon_D^1) \right\rceil + 1.
  \]
\end{lemma}
\begin{proof}%
  \Cref{cor:outer-iterates-well-defined} ensures that the iterates \((x_k, y_k)\) are well defined and satisfy~\eqref{eq:accuracy-condition}.
  From \Cref{alg:exactpen-outer-loop} and~\eqref{eq:accuracy-update}, \(\epsilon_P^k\) drops below \(\epsilon\) in at most
  \(\lceil \log(\epsilon/\epsilon_P^1) / \log(\mu_{\epsilon}) \rceil\)
  iterations.
  This implies that \(\tau_k\) remains constant for at most
  \[
    \lceil \log(\epsilon / \epsilon_P^1) / \log(\mu_{\epsilon}) \rceil +
    \lceil \log(\epsilon / \epsilon_D^1) / \log(\mu_{\epsilon}) \rceil
  \]
  outer loop iterations.
  Otherwise, we would have both \(\epsilon_P^k\) and \(\epsilon_D^k\) below \(\epsilon\) which implies~\eqref{eq:overall-stopping-criterion} from~\eqref{eq:accuracy-condition} and \(\|c(x_k)\|_{\infty} \leq \epsilon_P^k\) so that \Cref{alg:exactpen-outer-loop} stops.

  Therefore, since \(\epsilon_P^k\) is never reset whereas \(\epsilon_D^k\) is reset each time \(\tau_k\) increases, \(\tau_k \geq \bar \tau\) in at most
  \[
    \left\lceil \frac{\bar \tau - \tau_{\mathrm{init}}}{\tau_{\min}} \right\rceil +
    \lceil \log(\epsilon / \epsilon_P^1) / \log(\mu_{\epsilon}) \rceil +
    \left\lceil \frac{\bar \tau - \tau_{\mathrm{init}}}{\tau_{\min}} \right\rceil \lceil \log(\epsilon / \epsilon_D^1) / \log(\mu_{\epsilon}) \rceil
  \]
  outer loops.
  Once \(\tau_k\) attains this threshold, \(\|c(x_k)\|_{\infty} \leq \epsilon\).
  Therefore, besides the iterations that drive \(\epsilon_P^k\) below \(\epsilon\), which are already counted, there are at most
  \[
    \lceil \log(\epsilon / \epsilon_D^1) / \log(\mu_{\epsilon}) \rceil
  \]
  additional iterations before \(\epsilon_D^k \leq \epsilon\).
  Once this happens,~\eqref{eq:overall-stopping-criterion} follows from~\eqref{eq:accuracy-condition} at the next outer iteration.
\end{proof}

In \Cref{lemma:outer-loop-conv}, we have shown that \(k_\epsilon = \mathcal{O}(|\log(\epsilon)|)\).
The assumption that such a \(\bar \tau\) exists in \Cref{lemma:outer-loop-conv} is motivated by the exactness of the penalty function stated in \Cref{thm:exactness}.
In general however, we should only expect \Cref{alg:exactpen-outer-loop} to reach an approximate local solution of~\eqref{eq:infeasibility-problem} rather than an approximate feasible solution when \(\tau_k \geq \bar \tau\), which is the assumption used in \citep{diouane-gollier-orban-2026, cartis-gould-toint-2011}.
When this approximate local solution is infeasible for~\eqref{eq:nlp}, the infeasibility detection mechanism of \Cref{sec:infeasibility-detection} would trigger.

Combining the results of \Cref{sec:inner-loop-convergence,sec:outer-loop-convergence}, we can compute a bound on the number of inner and outer loop iterations of \Cref{alg:exactpen-outer-loop}.
Each inner loop iteration is called with a tolerance \(\epsilon_D^k\), therefore, this number is
\[
  \sum_{k = 1}^{k_\epsilon} |\mathcal{S}_k(\epsilon_D^k)| + |\mathcal{U}_k(\epsilon_D^k)| + 1.
\]
\begin{theorem}%
  \label{thm:total-complexity}
  Let the assumptions of \Cref{lemma:outer-loop-conv} hold.
  For each \(k \geq 1\), let \(x_{k,0}\) denote the starting point of \Cref{alg:exactpen-inner-loop} at outer iteration \(k\) of \Cref{alg:exactpen-outer-loop}, and assume that \(\{x_{k,0}\}\) is bounded.
  Then,
  \begin{equation}%
    \label{eq:total-iteration-bound}
    \sum_{k = 1}^{k_\epsilon} |\mathcal{S}_k(\epsilon_D^k)| + |\mathcal{U}_k(\epsilon_D^k)| + 1 =
    \mathcal{O}(|\log(\epsilon)| \epsilon^{-2}).
  \end{equation}
\end{theorem}
\begin{proof}%
  Since \(f\) and \(c\) are continuous and the set \(\{x_{k,0}\}\) is bounded, \({(f + \bar \tau \|c\|_2)}_{\mathrm{up}} \coloneqq \sup_{k \geq 1} f(x_{k,0}) + \bar \tau \|c(x_{k,0})\|_2\) is finite.
  Since \(\tau_k \in [\tau_1, \bar \tau]\), \(\epsilon_D^k \geq \epsilon\),~\eqref{eq:S-bound} implies
  \begin{equation}
    \label{eq:S-bound-eps}
    |\mathcal{S}_k(\epsilon_D^k)| \leq
    \frac{
    {(f + \bar \tau \|c\|_2)}_{\mathrm{up}} - {(f + \tau_1 \|c\|_2)}_{\mathrm{low}}
    }{
    \tfrac{1}{2} \eta_1 \eta_C
    }
    \epsilon^{-2}.
  \end{equation}
  The result follows by combining~\eqref{eq:S-bound-eps} with~\eqref{eq:U-bound} and \Cref{lemma:outer-loop-conv} and observing that the right-hand side of~\eqref{eq:S-bound-eps} is independent of \(k\).
\end{proof}

\section{Step Computation}%
\label{sec:step-computation}

Now, we present our approach for solving~\eqref{eq:quadratic-subproblem}.
For future reference, we define the implicit functions \(s_l(\alpha)\) and \(y_{l+1}(\alpha)\) with
\begin{equation}%
  \label{eq:alpha-system}
  \begin{pmatrix}
    H_l + \sigma_l I & { J_l }^T \\
    J_l              & -\alpha I
  \end{pmatrix}
  \begin{pmatrix}
    s_l(\alpha) \\
    y_{l+1}(\alpha)
  \end{pmatrix}
  = -
  \begin{pmatrix}
    \nabla f_l \\
    c_l
  \end{pmatrix},
  \quad
  \text{for} \ \alpha \geq 0.
\end{equation}
When~\eqref{eq:alpha-system} has multiple solutions, which may happen for \(\alpha = 0\) when \(J_l\) is rank deficient or for \(\alpha \geq 0\) when \(H_l + \sigma_l I\) is singular, \(y_{l+1}(\alpha)\) denotes the least-norm \(y\)-component among them.
If \(H_l + \sigma_l I \succ 0\), all solutions share the same \(s_l(\alpha)\).

We separate our step computation procedure into two cases.
We first assume that \(H_l + \sigma_l I\) is positive definite and \(J_l\) has full row rank.
We call this the \emph{definite case}.
Then, we describe the general case with both assumptions lifted.

\subsection{Definite Case}%
\label{sec:step-computation-positive}

The following result characterizes the solution of~\eqref{eq:quadratic-subproblem}.
It is an adaptation of \citep[Corollary 5.4]{diouane-gollier-orban-2026}, which works on the Schur complement of the linear system in~\eqref{eq:alpha-system}, whereas we work on the full system.
Moreover, \citep[Corollary 5.4]{diouane-gollier-orban-2026} uses the notation \(\nu_l \coloneqq 1/ \sigma_l\), multiplies the first row of~\eqref{eq:alpha-system} by \(\nu_l\), and uses \(\tilde y_{l+1}(\alpha) \coloneqq \nu_l y_{l+1}(\alpha)\), which gives the equivalent system
\[
  \begin{pmatrix}
    \nu_l H_l  + I & { J_l }^T        \\
    J_l            & -\tilde \alpha I
  \end{pmatrix}
  \begin{pmatrix}
    s_l(\tilde \alpha) \\
    \tilde y_{l+1}(\tilde \alpha)
  \end{pmatrix}
  = -
  \begin{pmatrix}
    \nu_l \nabla f_l \\
    c_l
  \end{pmatrix},
  \quad
  \text{for} \ \tilde \alpha \geq 0,
\]
where \(\tilde \alpha \coloneqq \alpha / \nu_l\).

\begin{proposition}[\protect{\citealp[Corollary 5.4]{diouane-gollier-orban-2026}}]%
  \label{proposition:step-characterization}
  Assume that \(H_l + \sigma_l I\) is positive definite.
  Then the unique solution \(s^*_l \coloneqq s_l(\alpha^*_l)\) of~\eqref{eq:quadratic-subproblem} satisfies~\eqref{eq:alpha-system},
  where \( \alpha^*_l \geq 0 \) is such that \( \alpha^*_l (\| y_{l+1}(\alpha^*_l) \|_2 - \tau_k) = 0 \) and \( \| y_{l+1}(\alpha^*_l) \|_2 \leq \tau_k \).
\end{proposition}

Our approach for computing steps is based on \Cref{proposition:step-characterization}.
We first compute \(s_l(0)\) and \(y_{l+1}(0)\) by solving~\eqref{eq:alpha-system} with \(\alpha = 0\).
Then, we check whether \(\|y_{l+1}(0)\|_2 \leq \tau_k\).
If this is not the case, we solve
\begin{equation}
  \label{eq:alpha-equation}
  \|y_{l+1}(\alpha)\|_2 = \tau_k,
\end{equation}
for \(\alpha > 0\).
To solve the equation, we apply Newton's method to the equivalent \emph{secular} equation \(1/\|y_{l+1}(\alpha)\|_2 - 1 / \tau_k = 0\), which turns out to be better suited for numerical computations.
This approach is motivated by the method of \citet{more-sorensen-1983}, and we refer to \citep[Section 5]{diouane-gollier-orban-2026} for more details.
For completeness, the following is an adaptation of \citep[Lemma 7.3.1]{conn-gould-toint-2000} that characterizes the iterates of the Newton method.
\begin{lemma}%
  \label{lemma:cgt-tr-more-sorensen-2}
  Let \(\bar \alpha \geq 0\) be such that the matrix in~\eqref{eq:alpha-system} is nonsingular.
  Assume that \(y_{l+1}(\bar \alpha) \neq 0\).
  Then, the secular function \(\xi : \alpha \to 1/\|y_{l+1}(\alpha)\|_2 - 1 / \tau_k\) is differentiable at \(\bar \alpha\).
  Moreover,
  \[
    \frac{\xi(\bar \alpha)}{\xi'(\bar \alpha)} =
    -\frac{
    (1 - \tau_k^{-1}\|y_{l+1}(\bar \alpha)\|_2) \|y_{l+1}(\bar \alpha)\|_2^2
    }{
    {y_{l+1}(\bar \alpha)}^T y'_{l+1}(\bar \alpha)
    },
  \]
  where \(y'_{l+1}(\bar \alpha)\) is the solution of~\eqref{eq:alpha-system} with right-hand side \((0, y_{l+1}(\bar \alpha))\).
\end{lemma}
In \Cref{lemma:cgt-tr-more-sorensen-2}, the assumption that~\eqref{eq:alpha-system} is nonsingular is satisfied for all \(\bar \alpha\) when \(H_l + \sigma_l I\) is positive definite and \(J_l\) has full row rank.
The convergence of Newton's method is guaranteed by the following lemma.
\begin{lemma}[\protect{\citealp[Lemma 7.3.2]{conn-gould-toint-2000}}]%
  \label{lemma:cgt-tr-more-sorensen}
  Assume that \(H_l + \sigma_l I\) is positive definite.
  Then, Newton's method applied to the secular equation \(\alpha \to 1/\|y_{l+1}(\alpha)\|_2 - 1 / \tau_k\) started from any \(\alpha \in (- \lambda_{\min}(J_l {(H_l + \sigma_l I)}^{-1} J_l^T), \alpha^*_l]\) produces a sequence of iterates that monotonically increase towards \(\alpha^*_l\),
  where \(\alpha^*_l\) is the unique nonnegative root to the equation.
\end{lemma}

Crucially, each iteration of Newton's method applied to the secular equation requires the solution of the symmetric indefinite linear system~\eqref{eq:alpha-system} for different values of \(\alpha\).
We use sparse \(LDL^T\) factorization algorithms to compute these.
In our implementation, the matrix factorizations are computed using the open-source solver MUMPS \citep{amestoy-duff-lexcellent-2000, amestoy-duff-lexcellent-koster-2001} or the factorization routine MA57 \citep{duff-2004} from the Harwell subroutine library.
Instead of requiring convergence to the exact solution of \(q(\cdot; x_l, \sigma_l)\), we truncate Newton's method when we find \(s_l \coloneqq s_l(\alpha_l)\) that satisfies~\eqref{eq:simple-decrease}, which significantly reduces the number of iterations of Newton's method in practice.

In practice, we find that we rarely need to perform an iteration of Newton's method.
That is, the choice \(\alpha_l = 0\) often satisfies~\eqref{eq:simple-decrease} so that we do not need to update it.
We summarize our step computation procedure when \(H_l + \sigma_l I \succ 0\) and \(J_l\) has full row rank in \Cref{alg:more-sorensen-step-definite}.

\begin{algorithm}[ht]%
  \caption{%
    \label{alg:more-sorensen-step-definite}
    Step Computation -- Definite Case.
  }
  \begin{algorithmic}[1]%

    \Require{}%
    Symmetric positive definite matrix \(H_l + \sigma_l I\);
    full row rank Jacobian \(J_l\);
    gradient \(\nabla f_l\);
    constraint values \(c_l\);
    penalty parameter \(\tau_k\).

    \State{}%
    Set \(\alpha_l \coloneqq 0\).

    \State{}%
    Solve
    \begin{equation}%
      \label{eq:kkt-system}
      \begin{pmatrix}
        H_l + \sigma_l I & { J_l }^T   \\
        J_l              & -\alpha_l I
      \end{pmatrix}
      \begin{pmatrix}
        s_l \\
        y_{l+1}
      \end{pmatrix}
      = -
      \begin{pmatrix}
        \nabla f_l \\
        c_l
      \end{pmatrix},
    \end{equation}

    \While{\(q(s_l; x_l, \sigma_l) > q(0; x_l, \sigma_l)\)}%

    \State{}%
    Solve
    \begin{equation}%
      \label{eq:kkt-system-derivative}
      \begin{pmatrix}
        H_l + \sigma_l I & { J_l }^T   \\
        J_l              & -\alpha_l I
      \end{pmatrix}
      \begin{pmatrix}
        s'_l \\
        y'_{l+1}
      \end{pmatrix}
      =
      \begin{pmatrix}
        0 \\
        y_{l+1}
      \end{pmatrix}.
    \end{equation}

    \State%
    Set
    \Comment{Newton step for the secular equation.}
    \begin{equation}%
      \label{eq:alpha-update}
      \alpha_{l}^+ \coloneqq \alpha_{l} + \frac{(1 - \tau_k^{-1}\| y_{l+1} \|_2) \| y_{l+1} \|_2^2}{{ y_{l+1} }^T y'_{l+1}}.
    \end{equation}

    \State%
    Set \(\alpha_{l} \coloneqq \alpha_{l}^+\) and solve~\eqref{eq:kkt-system}.

    \EndWhile%

  \end{algorithmic}
\end{algorithm}

In the definite case, the step computed by \Cref{alg:more-sorensen-step-definite} satisfies \Cref{stepassumption:step-s-stationarity,stepassumption:decrease}.
Condition~\eqref{eq:simple-decrease} is satisfied by construction, and~\eqref{eq:cauchy-decrease} follows from~\eqref{eq:simple-decrease} by \Cref{lemma:simple-cauchy-decrease}, so that \Cref{stepassumption:decrease} holds.
Moreover, since \((s_l, y_{l+1})\) solves~\eqref{eq:kkt-system}, its first block row is~\eqref{eq:step-s-stationarity}, which gives \Cref{stepassumption:step-s-stationarity}.

\subsection{General Case}%
\label{sec:step-computation-negative}

In general, the positive definiteness of \(H_l+\sigma_l I\) and the full-rank condition on \(J_l\) may not hold.
To detect potential degeneracy, we use the inertia of~\eqref{eq:kkt-system} as an indicator for the properties of both matrices.
Specifically, the inertia is the triple \((n_+,n_-,n_0)\), where \(n_+\), \(n_-\), and \(n_0\) denote the numbers of positive, negative, and zero eigenvalues, respectively.
When a sparse factorization routine is used to solve~\eqref{eq:kkt-system}, the inertia is a by-product of the factorization.

\begin{proposition}%
  \label{prop:inertia}
  Assume that \(H_l + \sigma_l I\) is positive definite.
  Then,
  \begin{enumerate}
    \item The inertia of~\eqref{eq:kkt-system} with \(\alpha_l = 0\) is \((n, r, m - r)\), where \(r = \textup{Rank}(J_l)\).
    \item The inertia of~\eqref{eq:kkt-system} with \(\alpha_l > 0\) is \((n, m, 0)\).
  \end{enumerate}
\end{proposition}
\begin{proof}
  This is a direct consequence of the Haynsworth inertia additivity formula \citep{haynsworth-1968}.
\end{proof}

Our approach for detecting degeneracy in \(H_l + \sigma_l I\) and \(J_l\) is based on \Cref{prop:inertia}.
First, we try to solve~\eqref{eq:kkt-system} with \(\alpha_l = 0\) and compute the inertia \((n_+, n_-, n_0)\).
If \(n_- < m\), then we set \(\alpha_l\) to a small positive value \(\alpha_{\textup{trial}} > 0\) and retry solving~\eqref{eq:kkt-system}, which updates \((n_+, n_-, n_0)\).
After this first correction step, if \(n_+ < n\), we increase \(\sigma_l\) and resolve~\eqref{eq:kkt-system} with the updated value of \(\alpha_l\) until \((n_+, n_-, n_0) = (n, m, 0)\).

Following this inertia correction step, two cases may occur: either \(H_l + \sigma_l I\) is positive definite or not.

First, if \(H_l + \sigma_l I\) is positive definite, the same procedure as \Cref{alg:more-sorensen-step-definite} applies to compute a step satisfying \Cref{stepassumption:step-s-stationarity,stepassumption:decrease}, regardless of rank deficiency in \(J_l\), provided that \(\alpha_{\textup{trial}} > 0\).
\begin{proposition}%
  \label{prop:step-computation-rank-deficient}
  Assume that \(H_l + \sigma_l I\) is positive definite and that \(J_l\) is rank deficient.
  Assume that \(s = 0\) is not a solution to~\eqref{eq:quadratic-subproblem}.
  Then, there is an \(\alpha_{\textup{trial}} > 0\) such that \Cref{alg:more-sorensen-step-definite} started with \(\alpha_l = \alpha_{\textup{trial}}\) computes a step \(s_l\) satisfying \Cref{stepassumption:step-s-stationarity,stepassumption:decrease}.
\end{proposition}
\begin{proof}
  \Cref{stepassumption:step-s-stationarity} is satisfied by construction, since \(s_l\) is the solution of~\eqref{eq:alpha-system} for some \(\alpha \geq 0\).
  Moreover, following \Cref{lemma:simple-cauchy-decrease}, it is sufficient to show \Cref{alg:more-sorensen-step-definite} produces a step \(s_l\) satisfying~\eqref{eq:simple-decrease}.

  Since \(H_l + \sigma_l I\) is positive definite, \Cref{proposition:step-characterization} applies unchanged.
  In \Cref{proposition:step-characterization}, either \(\alpha^*_l = 0\) or \(\alpha^*_l > 0\).

  If \(\alpha^*_l = 0\), then the solution of~\eqref{eq:quadratic-subproblem} is \(s^*_l = s_l(0)\).
  In particular, since \(s_l^* \neq 0\) and \(s_l^*\) is the unique solution of~\eqref{eq:quadratic-subproblem}, by strict convexity of~\eqref{eq:quadratic-subproblem},
  \[
    q(s_l(0); x_l, \sigma_l) < q(0; x_l, \sigma_l).
  \]
  Continuity of \(\alpha \to q(s(\alpha); x_l, \sigma_l)\) on \([0, \infty)\) then implies that there is an \(\alpha_{\textup{trial}} > 0\) such that \(s_l\) satisfies~\eqref{eq:simple-decrease} for all \(\alpha_l \in [0, \alpha_{\textup{trial}}]\).

      If \(\alpha^*_l > 0\), then if \(\alpha_{\textup{trial}} \in (0, \alpha^*_l]\), \Cref{lemma:cgt-tr-more-sorensen} implies that the Newton iterates~\eqref{eq:alpha-update} converge to \(\alpha^*_l\).
  In particular, \Cref{alg:more-sorensen-step-definite} converges to the solution of~\eqref{eq:quadratic-subproblem} and produces an \(s_l\) satisfying~\eqref{eq:simple-decrease}.
\end{proof}
The assumption that \(s = 0\) is not a solution to~\eqref{eq:quadratic-subproblem} in \Cref{prop:step-computation-rank-deficient} is not restrictive.
Indeed, if the solution of~\eqref{eq:quadratic-subproblem} is \(s = 0\), following \Cref{proposition:step-characterization}, there is an \(\alpha_l^* \geq 0\) such that \((0, y_{l+1}(\alpha_l^*))\) is a solution of~\eqref{eq:alpha-system}.
The first block row of~\eqref{eq:alpha-system} then implies that \(\nabla f_l + J_l^T y_{l+1}(\alpha_l^*) = 0\), which means that \(x_l\) is stationary with respect to~\eqref{eq:kkt-dual}, in which case, \Cref{alg:exactpen-inner-loop} terminates and returns \(x_l\).

If after the inertia correction step, \(H_l + \sigma_l I\) is not positive definite, a possible approach is to increase \(\sigma_l\) and restart the procedure until it is.
However, experience showed that this usually leads to unnecessarily large values in \(\sigma_l\) and slower convergence of the algorithm.
Indeed, observe that if \(\alpha_l = 0\), then even if \(H_l + \sigma_l I\) is indefinite, when the inertia of~\eqref{eq:alpha-system} is \((n, m, 0)\), the solution \(s_l\) of~\eqref{eq:alpha-system} is a solution of
\[
  \min_s f_l + {\nabla f_l}^T s + \tfrac{1}{2} s^T (H_l + \sigma_l I) s
  \quad \text{s.t.} \quad c_l + J_l s =  0,
\]
which often turns out to be a good step for~\eqref{eq:penalty-nlp} in practice.
We found that a better technique is to perform a very limited number \(j_{\max}\) of iterations of Newton's method using the update rule of \Cref{lemma:cgt-tr-more-sorensen-2}.
We stress that neither \Cref{proposition:step-characterization} nor \Cref{lemma:cgt-tr-more-sorensen} apply.
Hence, we cannot guarantee that we will satisfy~\eqref{eq:simple-decrease}.
Moreover, \Cref{lemma:simple-cauchy-decrease} no longer applies, so that~\eqref{eq:simple-decrease} does not imply~\eqref{eq:cauchy-decrease}.
As a consequence, we verify both~\eqref{eq:simple-decrease} and~\eqref{eq:cauchy-decrease}.
If both conditions are satisfied, we return the step; otherwise, we either continue to the next Newton iteration or increase \(\sigma_l\) and restart the procedure.
Note that \Cref{stepassumption:step-s-stationarity} is satisfied by construction, since \(s_l\) is the solution of~\eqref{eq:alpha-system}.

\begin{algorithm}[ht]%
  \caption{%
    \label{alg:more-sorensen-step-indefinite}
    Step Computation -- General Case.
  }
  \begin{algorithmic}[1]%

    \Require{}%
    Symmetric matrix \(H_l\);
    regularization parameter \(\sigma_l > 0\);
    Jacobian \(J_l\);
    gradient \(\nabla f_l\);
    constraint values \(c_l\);
    penalty parameter \(\tau_k\);
    constants \(\alpha_{\mathrm{trial}}^1 > 0\);
    \(\alpha_{\mathrm{trial}}^2 > \alpha_{\mathrm{trial}}^1\);
    \(\mu_{\sigma} > 1\);
    \(\mu_{\alpha} \in (0, 1)\);
    \(\sigma_{\mathrm{ub}} > 0\);
    \(j_{\max} \in \mathbb{N}\);
    \(\eta_C > 0\).

    \State{}%
    Set \(\alpha_l \coloneqq 0\).

    \State{}%
    \label{algline:general-solve}
    Attempt to solve~\eqref{eq:kkt-system}, compute the inertia \((n_+, n_-, n_0) \).
    \Comment{Inertia corrections.}

    \State{}%
    If either the factorization~\eqref{eq:kkt-system} fails or \(n_- < m\), set \(\alpha_l \coloneqq \alpha_{\mathrm{trial}}^1\) and retry solving~\eqref{eq:kkt-system}.
    If the factorization fails again or if \(n_- < m\), set \(\alpha_l \coloneqq \alpha_{\mathrm{trial}}^2\) and retry solving~\eqref{eq:kkt-system}.
    If it still fails, proceed to \Cref{algline:general-sigma-increase}.

    \State{}%
    \label{algline:general-sigma-increase}
    Until \(n_+ = n\) and the factorization~\eqref{eq:kkt-system} succeeds, set \(\sigma_l \coloneqq \mu_{\sigma}\sigma_{l}\), repeat \Cref{algline:general-solve}.
    If \(\sigma_l > \sigma_{\mathrm{ub}}\), \textbf{set} \(H_l \coloneqq 0\) and \textbf{restart}.

    \State{}%
    Compute the inertia of \(H_l + \sigma_l I\) with an \(LDL^T\) factorization.

    \State{}%
    If \(H_l + \sigma_l I\) is positive definite, compute \(s_l\) with \Cref{alg:more-sorensen-step-definite} and \textbf{return}.

    \For{\(j = 1, \dots, j_{\max}\)}%

    \State{}%
    \label{algline:general-return}
    If \(s_l\) satisfies~\eqref{eq:simple-decrease} and~\eqref{eq:cauchy-decrease}, \textbf{return}.

    \State{}%
    Solve~\eqref{eq:kkt-system-derivative}.

    \State{}%
    Compute \(\alpha_{l}^+\) with~\eqref{eq:alpha-update}.

    \State{}%
    \label{algline:general-alpha-update}
    If \(\alpha_{l}^+ < 0\), set \(\alpha_{l} \coloneqq \max(\alpha_{\mathrm{trial}}^1, \mu_{\alpha} \alpha_l)\), else set \(\alpha_{l} \coloneqq \alpha_{l}^+\). Solve~\eqref{eq:kkt-system}.

    \EndFor%

    \State{}%
    Set \(\sigma_l \coloneqq \mu_{\sigma}\sigma_{l}\) and \textbf{restart}.

  \end{algorithmic}
\end{algorithm}

In practice, before committing to the \(LDL^T\) factorization of \(H_l + \sigma_l I\), we try cheap tests to check whether \(H_l + \sigma_l I\) is positive definite.
Specifically, we perform the following two tests:
\begin{enumerate}
  \item We verify whether all diagonal elements of \(H_l + \sigma_l I\) are positive. We return that \(H_l + \sigma_l I\) is not positive definite if this is not the case.
  \item We verify whether \(H_l + \sigma_l I\) is strictly diagonally dominant. We return that \(H_l + \sigma_l I\) is positive definite if this is the case.
\end{enumerate}
This allows us to save a few \(LDL^T\) factorizations in practice.

Note that increasing \(\sigma_l\) in \Cref{alg:more-sorensen-step-indefinite} does not break the convergence of \Cref{alg:exactpen-inner-loop}.
Indeed, the increase of \(\sigma_l\) in \Cref{alg:more-sorensen-step-indefinite} is proportional to a constant factor \(\mu_{\sigma} > 1\), which is the same mechanism as in \Cref{alg:exactpen-inner-loop}.

Finally, the value \(\alpha_k\) used in the extrapolation technique of \Cref{sec:extrapolation-technique} is the last value of \(\alpha_l\) computed by \Cref{alg:more-sorensen-step-definite} or \Cref{alg:more-sorensen-step-indefinite} in the final iteration of \Cref{alg:exactpen-inner-loop} at outer iteration \(k\).

\section{Implementation Details}%
\label{sec:implementation-details}

In addition to the core algorithmic components presented in the previous sections, our implementation incorporates several practical strategies to improve performance and robustness.
This section describes and justifies these heuristics, which include nonmonotone strategies, initialization techniques, and regularization parameter updates.

\subsection{Nonmonotone Strategies}%
\label{sec:nonmonotone-strategies}

A first practical issue encountered in our implementation arises when the penalty parameter \(\tau_k\) significantly exceeds the threshold \(\bar \tau\) presented in \Cref{thm:exactness}.
For example, this can occur when \(\tau_1 = \max(\tau_{\mathrm{init}}, \|y_0\|_1) \gg \bar \tau\).
In view of the right-hand side of~\eqref{eq:S-bound}, the number of inner loop iterations increases with \(\tau_k\).
Moreover, similar to the observations of \citet{maratos-1978}, we find that the algorithm requires smaller steps, that is, larger values of \(\sigma_l\), when \(\tau_k\) is too large.
To mitigate this effect, we implement two concurrent nonmonotone strategies that relax the step acceptance rule~\eqref{eq:armijo} in \Cref{alg:exactpen-inner-loop}.

The first strategy, proposed by \citet{grippo-lampariello-lucidi-1989}, modifies the value of \(\rho_l\) in~\eqref{eq:armijo}.
Let \(q \in \N\) be a user-defined memory parameter.
Let \(q_l \coloneqq \min(q, l)\) if \(q > 0\) and \(q_l \coloneqq 0\) otherwise.
Define \(\mathcal{S}^+_{q_l}\) as the set of the \(q_l\) most recent successful iterations.
We define
\begin{equation}%
  \label{eq:nonmonotone-max}
  {(f + \tau_k \|c\|_2)}_{\max, l} \coloneqq \max_{j \in \mathcal{S}^+_{q_l}} f(x_j) + \tau_k \|c(x_j)\|_2.
\end{equation}
For completeness, we set \(\mathcal{S}^+_0 = \{ 0 \}\), and define \({(f + \tau_k \|c\|_2)}_{\max, 0} \coloneqq f(x_0) + \tau_k \|c(x_0)\|_2\).
Then, we modify~\eqref{eq:armijo} with
\begin{equation}%
  \label{eq:armijo-nonmonotone}
  \rho_l \coloneqq
  \frac{
  {(f + \tau_k \|c\|_2)}_{\max, l} - \left( f(x_l + s_l) + \tau_k \|c(x_l + s_l)\|_2 \right)
  }{
  {(f + \tau_k \|c\|_2)}_{\max, l} - (f_l + {\nabla f_l}^T s_l + \tau_k\|c_l + J_l s_l\|_2)
  } \geq \eta_1.
\end{equation}
To disable this heuristic, it suffices to set \(q = 0\).

While the strategy above has proved very effective in practice, we still observed a few instances where \Cref{alg:exactpen-inner-loop} struggles when the penalty parameter is too large.
To address this, we implemented a second nonmonotone strategy.
Our approach is inspired by the watchdog technique proposed by \citet{chamberlain-lemarechal-pedersen-powell-1982}.
Specifically, it allows a few consecutive iterations to be accepted even if they do not satisfy~\eqref{eq:armijo-nonmonotone}.

The watchdog has two states: \textit{active} and \textit{inactive}.
It is \textit{inactive} when \Cref{alg:exactpen-inner-loop} starts.
Suppose that, at some iterate \(x_{l_0}\), the watchdog is \textit{inactive} and the step \(s_{l_0}\) satisfies \(\rho_{l_0} < \eta_1\) in~\eqref{eq:armijo-nonmonotone}.
Instead of rejecting the step, we check whether it sufficiently decreases either \(f\) or the penalty term, taken separately.
More specifically, we compute
\[
  \rho_{l_0}^f \coloneqq
  \frac{
  f(x_{l_0}) - f(x_{l_0} + s_{l_0})
  }{
  - {\nabla f(x_{l_0})}^T s_{l_0}
  }
  \quad \text{and} \quad
  \rho_{l_0}^c \coloneqq
  \frac{
    \|c(x_{l_0})\|_2 - \|c(x_{l_0} + s_{l_0})\|_2
  }{
    \|c(x_{l_0})\|_2 - \|c(x_{l_0}) + J(x_{l_0}) s_{l_0}\|_2
  },
\]
where we set \(\rho_{l_0}^f \coloneqq 0\) if \({\nabla f(x_{l_0})}^T s_{l_0} > 0\), and \(\rho_{l_0}^c \coloneqq 0\) if \(\|c(x_{l_0}) + J(x_{l_0}) s_{l_0}\|_2 > \|c(x_{l_0})\|_2\), since the corresponding linear model then predicts no decrease.
If \(\max(\rho_{l_0}^f, \rho_{l_0}^c) \geq \eta_3\), where \(\eta_3 \in (0, 1)\) is a user-defined parameter, we activate the watchdog, accept the step \(s_{l_0}\), keep \(\sigma_{l_0 + 1} = \sigma_{l_0}\) and continue with \Cref{alg:exactpen-inner-loop}.
Otherwise, we reject the step and proceed to the next iteration of \Cref{alg:exactpen-inner-loop}.
While the watchdog is active, we accept steps based on the monotone ratio~\eqref{eq:armijo} instead of~\eqref{eq:armijo-nonmonotone}.

Once activated, the watchdog either succeeds or fails.
If it remains active for more than \(r \in \N\) iterations, it fails.
In this case, we deactivate it, reset \((x_l, y_l)\) to \((x_{l_0}, y_{l_0})\) and increase \(\sigma_l \in [\gamma_1^r\sigma_{l_0}, \gamma_2^r\sigma_{l_0}]\).
On the other hand, if for some \(r_{l_0} \leq r\) and \(l \coloneqq l_0 + r_{l_0}\), the step \(s_{\mathrm{wd}, l} \coloneqq x_l - x_{l_0}\) satisfies either
\begin{equation}%
  \label{eq:watchdog-success-1}
  f(x_{l_0}) + \tau_k \|c(x_{l_0})\|_2 - f(x_{l_0} + s_{\mathrm{wd}, l}) - \tau_k \|c(x_{l_0} + s_{\mathrm{wd}, l})\|_2 \geq \tfrac{1}{2}\eta_0 \min(1, \|\nabla f_l + J_l^T y_{l+1}\|_\infty^2),
\end{equation}
or
\begin{equation}%
  \label{eq:watchdog-success-2}
  \|\nabla f_l + J_l^T y_{l+1} \|_{\infty} \leq (1 - \eta_0) \|\nabla f_{l_0} + J_{l_0}^T y_{l_0+1}\|_{\infty},
\end{equation}
where \(\eta_0 \in (0, 1)\) is a user-defined parameter, then it succeeds and is deactivated.
We justify~\eqref{eq:watchdog-success-1} with \Cref{stepassumption:decrease}.
For~\eqref{eq:watchdog-success-2}, a uniform decrease in the optimality measure \(\|\nabla f_l + J_l^T y_{l+1} \|_{\infty}\) should eventually lead to \(\|\nabla f_l + J_l^T y_{l+1} \|_{\infty} \to 0\).

\subsection{Quadratic Regularization Parameter Heuristics}%
\label{sec:regularization-updates}

As explained in \Cref{sec:penalized-subproblem}, and shown in \Cref{sec:inner-loop-convergence}, the choice of the regularization parameter~\(\sigma_l\) is crucial for both the practical performance and the theoretical convergence of \Cref{alg:exactpen-inner-loop}.
We now discuss how we update \(\sigma_l\) throughout the inner iterations.

When we first call \Cref{alg:exactpen-inner-loop}, on the first outer iteration, we initialize the regularization parameter as
\begin{equation}%
  \label{eq:sigma0-init}
  \sigma_0 = \epsilonm.
\end{equation}
In subsequent calls to \Cref{alg:exactpen-inner-loop}, we initialize \(\sigma_0\) at outer iteration \(k + 1\) as the final value of \(\sigma_l\) on outer iteration~\(k\).
To ensure that the convergence analysis in \Cref{sec:inner-loop-convergence} remains applicable, we project the former into the interval \([\sigma_{\min}, \sigma_{\mathrm{ub}}]\) with \(\sigma_{\min} = \epsilonm\) and the same value of \(\sigma_{\mathrm{ub}}\) as in \Cref{alg:more-sorensen-step-indefinite}.
In particular, the value of \(\sigma_{\max}\) in~\eqref{eq:sigma-max} can be bounded independently of \(k\).
Then, \(\sigma_{0}\) in the upper bound of \Cref{lemma:U-complexity} can be replaced with \(\sigma_{\min}\) which is independent of \(k\) too.

The motivation for~\eqref{eq:sigma0-init} is to be aggressive near a potential solution.
Indeed, we have observed fast local convergence behavior when \(\sigma_l \approx 0\).
Therefore, when the initial pair \((x_0, y_0)\) is close to a solution of~\eqref{eq:nlp}, the choice~\eqref{eq:sigma0-init} allows us to converge to high levels of accuracy in just a few iterations.

When this is not the case however, we often observe that the aggressive initialization~\eqref{eq:sigma0-init} produces poor agreement between the actual and first-order predicted reductions measured in~\eqref{eq:armijo-nonmonotone}.
As a result, the ratio \(\rho_l\) in~\eqref{eq:armijo-nonmonotone} is often negative during the first few inner iterations.
Applying the standard update strategy in this situation generates a sequence of unsuccessful iterations to reach higher orders of magnitude for \(\sigma_l\).
To avoid this behavior, our implementation allows for a larger increase of \(\sigma_l\) on the first unsuccessful iteration of \Cref{alg:exactpen-inner-loop}.
That is, when we encounter the first unsuccessful iteration, we set
\(
\sigma_{l+1} \coloneqq \max(\gamma_2\sigma_{l}, \sqrt{\|\nabla f_l + J_l^T y_{l+1}\|_{\infty}}).
\)
Choosing \(\sigma_{l+1}\) proportional to the square root of a measure of dual infeasibility has recently been suggested in \citep{mishchenko-2023,doikov-nesterov-2024}.
We observed that this choice saves a few unsuccessful iterations during the early stages of the algorithm.
Moreover, this strategy does not affect the proof mechanism developed in \Cref{sec:inner-loop-convergence} because it only allows for a larger increase of \(\sigma_l\) once.

\subsection{Linear Algebra Routines}%
\label{sec:linear-algebra-routines}

By default, we solve the linear systems~\eqref{eq:kkt-system} using the multifrontal symmetric-indefinite \(LDL^T\) factorization implemented in MUMPS~\citep{amestoy-duff-lexcellent-2000,amestoy-duff-lexcellent-koster-2001,montoison-orban-sweeney-2024}.
We use a default pivot tolerance of
\(
\epsilon_{\mathrm{piv}} = \epsilonm,
\)
which favors sparse factors over pivoting robustness and results in faster factorizations.

We also use iterative refinement, with at most \(5\) refinement steps.
After each solve, we compare the backward error estimate reported by MUMPS with \(\sqrt{\epsilonm}\).
If the estimate exceeds this threshold, we consider the solve unsuccessful, increase the pivot tolerance to \(10^{-2}\), and force exactly \(5\) refinement steps for all subsequent solves.

This adaptive strategy exploits the fact that well-conditioned systems can typically be solved efficiently with a small pivot tolerance, while a single ill-conditioned system warrants more conservative factorization parameters.
To avoid using an unnecessarily large pivot tolerance indefinitely, we reduce \(\epsilon_{\mathrm{piv}}\) by a factor of \(10\) whenever the backward error estimate falls below \(100\epsilonm\) for more than \(20\) consecutive solves.

Finally, we use MUMPS automatic scaling of~\eqref{eq:kkt-system} to improve the conditioning of the linear systems.

\subsection{Initialization}%
\label{sec:initialization}

For the initial approximation of the Lagrange multiplier \(y_0\), if no initial guess is provided, we initialize \(y_0\) as the least-squares solution of \(J_0^T y = - \nabla f_0\) by solving
\begin{equation}%
  \label{eq:least-squares-init}
  \begin{pmatrix}
    I   & J_0^T \\
    J_0 & 0
  \end{pmatrix}
  \begin{pmatrix}
    w \\
    y_0
  \end{pmatrix}
  =
  \begin{pmatrix}
    - \nabla f_0 \\
    0
  \end{pmatrix}.
\end{equation}
We use the same sparse \(LDL^T\) factorization routines as in \Cref{sec:linear-algebra-routines} to solve~\eqref{eq:least-squares-init}.
If the factorization of the matrix in~\eqref{eq:least-squares-init} fails, we retry by setting the lower-right block to \(-\alpha_{\textup{trial}}^1 I\) first and then \(-\alpha_{\textup{trial}}^2 I\) as in \Cref{alg:more-sorensen-step-indefinite}.

If \(\nabla f_0 = 0\), however, \(y_0 = 0\), so that, in view of \Cref{line:init-tau} of \Cref{alg:exactpen-outer-loop}, the value of \(\tau_1\) is set to \(\tau_{\mathrm{init}}\).
For example, this happens when \(f\) is constant.
In this case, instead of solving~\eqref{eq:least-squares-init}, we solve
\begin{equation}%
  \label{eq:least-squares-init-2}
  \begin{pmatrix}
    I   & J_0^T \\
    J_0 & 0
  \end{pmatrix}
  \begin{pmatrix}
    w \\
    y_0
  \end{pmatrix}
  =
  \begin{pmatrix}
    0 \\
    - c_0
  \end{pmatrix}.
\end{equation}

The advantage of~\eqref{eq:least-squares-init} and~\eqref{eq:least-squares-init-2} is that if \(x_0\) is a local solution of~\eqref{eq:nlp}, then our algorithm stops immediately after having computed \(y_0\).
However, a problem with the initialization strategy above is that it may produce a very large value of \(\|y_0\|_1\), especially when \(J_0\) is ill-conditioned, which in turn produces a very large value of \(\tau_1\), in view of \Cref{line:init-tau} of \Cref{alg:exactpen-outer-loop}.
To avoid this, similarly to IPOPT~\citep[Section 3.6]{wachter-biegler-2006}, we first implemented a simple safeguard that checks whether \(\|y_0\|_1 > y_{\max}\) and discards \(y_0\) when this condition is met.
This strategy, however, only improved the performance of the algorithm on a few problems, while unnecessarily discarding \(y_0\) on many others.
We found that a numerically more robust approach is to detect ill-conditioning of either~\eqref{eq:least-squares-init} or~\eqref{eq:least-squares-init-2}.
Let \(K\) be the matrix in~\eqref{eq:least-squares-init} or~\eqref{eq:least-squares-init-2}, \(u\) the right-hand side of the corresponding system and \(\|K\|_{\max}\) the maximum absolute value of the entries of \(K\).
If \(K\) is nonsingular, then its condition number satisfies
\[
  \kappa_{\infty}(K) = \|K\|_\infty \|K^{-1}\|_\infty \geq \|K\|_{\max} \|(w, y_0)\|_{\infty} / \|u\|_{\infty}.
\]
We therefore check whether
\begin{equation}%
  \label{eq:condition-number}
  \|K\|_{\max} \|(w, y_0)\|_{\infty} / \|u\|_{\infty} > 1 / \epsilonm.
\end{equation}
When~\eqref{eq:condition-number} is satisfied, we retry the initialization by setting the lower-right block of \(K\) to \(-\alpha_{\textup{trial}}^1 I\) first and then \(-\alpha_{\textup{trial}}^2 I\) as in \Cref{alg:more-sorensen-step-indefinite}.
If~\eqref{eq:condition-number} is met for both \(\alpha_{\textup{trial}}^1\) and \(\alpha_{\textup{trial}}^2\), we set \(y_0 = 0\).

\subsection{Quasi-Newton Approximations}%
\label{sec:quasi-newton-approximations}

For the definition of \(H_l\) in~\eqref{eq:quadratic-subproblem}, our algorithm allows for quasi-Newton approximations of the Hessian of the Lagrangian.
Namely, our algorithm implements a limited memory BFGS method to approximate \(H_l \approx \nabla^2 \mathcal{L} (x_l, y_l)\).
At each inner iteration \(l\), we compute
\[
  y_{\mathrm{qn}, l} \coloneqq \nabla \mathcal{L}(x_{l+1}, y_{l+1}) - \nabla \mathcal{L}(x_l, y_{l + 1}).
\]
When we accept the step \(s_l\), we push the pair \((s_l, y_{\mathrm{qn}, l})\) into a limited memory BFGS data structure using the compact representation of \citet{byrd-nocedal-schnabel-1994}.
This allows us to write \(H_l\) as
\[
  H_l = \delta_l I + N_l M_l N_l^T.
\]
Our implementation defines \(\delta_l\) as \(\delta_l = y_{\mathrm{qn}, l}^T y_{\mathrm{qn}, l} / s_l^T y_{\mathrm{qn}, l}\).
To ensure positive definiteness and boundedness of \(H_l\), the update is skipped unless
\[
  s_l^T y_{\mathrm{qn}, l} > \epsilonm \|y_{\mathrm{qn}, l}\|_2 \max(\|s_l\|_2, \|y_{\mathrm{qn}, l}\|_2).
\]
The condition \(s_l^T y_{\mathrm{qn}, l} > \epsilonm \|s_l\|_2 \|y_{\mathrm{qn}, l}\|_2\) ensures positive definiteness, while the condition \(s_l^T y_{\mathrm{qn}, l} > \epsilonm \|y_{\mathrm{qn}, l}\|_2^2\) ensures that \(H_l\) is bounded.
If more than two consecutive updates are skipped, the approximation is reset.

In \Cref{alg:more-sorensen-step-definite,alg:more-sorensen-step-indefinite}, when we solve~\eqref{eq:kkt-system} and~\eqref{eq:kkt-system-derivative}, we use a standard Sherman-Morrison-Woodbury formula.
That is, we use an \(LDL^T\) factorization of the \(2 \times 2\) block system with \(H_l\) replaced by \(\delta_l I\) and then update the inverse by adding the low-rank term \(N_l M_l N_l^T\), following the procedure described in \citep[Section 3.8]{waltz-morales-nocedal-orban-2006}.

\subsection{Handling Small Steps}%
\label{sec:small-search-directions}

In a few cases, we have observed that the steps \(s_l\) computed by \Cref{alg:more-sorensen-step-indefinite} become very small relative to \(x_l\).
For this reason, at each inner iteration \(l\), we check whether
\begin{equation}%
  \label{eq:small-step}
  | s_l^{(i)} | < \epsilonm (1 + | x_l^{(i)} |), \quad \textup{for each} \ i = 1, \dots, n.
\end{equation}
When this happens, the solver effectively stalls, therefore we terminate the overall algorithm with a corresponding failure message.
To make sure that~\eqref{eq:small-step} is not triggered by some numerical error, we return the failure message when~\eqref{eq:small-step} is satisfied for \(n_{\mathrm{small}} + 1\) consecutive inner iterations.

\subsection{Other Heuristics}%
\label{sec:other-heuristics}

By default, our implementation is not invariant to scaling of the variables, objective function, and constraints.
Similar behavior has been reported for other solvers, including IPOPT~\citep[Section 3.8]{wachter-biegler-2006}.
We implemented the scaling strategy described in~\citep[Section 3.8]{wachter-biegler-2006}, but found that it degraded performance of the method, therefore we did not include it in the implementation.

We use two safeguards to adapt the penalty parameter \(\tau_k\).
First, for some values of \(\tau_k\),
\(
x \to f(x) + \tau_k \|c(x)\|_2
\)
can be unbounded below.
In this case, the iterates of \Cref{alg:exactpen-inner-loop} may satisfy \(f_l \to -\infty\).
We therefore terminate the inner loop with a failure message when
\(
f_l \leq -1/\epsilonm \approx -10^{16}.
\)
The outer loop then increases the penalty parameter according to
\(
\tau_{k+1} = 10\tau_k
\)
and restarts the inner loop.

Second, \(\tau_k\) can occasionally be substantially smaller than the penalty threshold discussed in \Cref{thm:exactness}.
In this situation, the algorithm may fail to sufficiently prioritize constraint feasibility, as indicated by
\(
\|c(x_k)\|_2 \geq \|c(x_{k-1})\|_2.
\)
When this occurs, we increase the penalty parameter more aggressively by multiplying the minimum penalty increase \(\tau_{\min}\) by a factor of \(10\).
When the penalty parameter becomes larger than a fixed threshold \(\tau_{\max} > 0\), we terminate the algorithm with a failure message.

Finally, to compute the tolerances \(\epsilon_D^k\) and \(\epsilon_P^k\) used in \Cref{alg:exactpen-outer-loop} in practice, instead of~\eqref{eq:accuracy-update} we use a relative tolerance based on the current residuals
\[
  \epsilon_D^k \coloneqq \min\left(1, \max\left(\mu_{\epsilon} \|\nabla f(x_{k-1}) + {J(x_{k-1})}^T y_{k-1}\|_\infty, \epsilon_D\right)\right) \quad \text{and} \quad \epsilon_P^k \coloneqq \min\left(1, \max\left(\mu_{\epsilon} \|c(x_{k-1})\|_\infty, \epsilon_P\right)\right).
\]
When either dual or primal infeasibility are already small, this avoids unnecessary calls to \Cref{alg:exactpen-inner-loop} with a loose tolerance that is satisfied immediately.

\subsection{Algorithmic Hyperparameters}%
\label{sec:algorithmic-hyperparameters}

This section summarizes all default algorithmic hyperparameters that we presented throughout this paper.
The default values of the constants in \Cref{alg:exactpen-outer-loop,alg:exactpen-inner-loop} are
\(\tau_{\mathrm{init}} = 1\);
\(\epsilon_P^1 = 1\);
\(\epsilon_D^1 = 1\);
\(\tau_{\min} = 1\);
\(\tau_{\max} = 1 / \epsilonm^{1.5} \approx 10^{23}\);
\(\mu_{\epsilon} = 0.01\);
\(\sigma_0 = \epsilonm \approx 10^{-16}\);
\(\eta_1 = \epsilonm^{0.25} \approx 10^{-4}\);
\(\eta_2 = 0.1\) by default and \(0.9\) when \(H_l\) is a quasi-Newton approximation; \(\gamma_1 = \gamma_2 = 3\);
\(\gamma_3 = 1/3\).
For the infeasibility detection described in \Cref{sec:infeasibility-detection}, we choose \(\epsilon_I = 10^{-3}\).
In \Cref{alg:more-sorensen-step-indefinite}, the default values of the constants are
\(\alpha_{\mathrm{trial}}^1 = \epsilonm^{0.8} \approx 10^{-13}\);
\(\alpha_{\mathrm{trial}}^2 = \epsilonm^{0.6} \approx 10^{-10}\);
\(\mu_{\alpha} = 0.1\);
\(\mu_{\sigma} = 10\);
\(\sigma_{\mathrm{ub}} = 1/\epsilonm^{0.8} \approx 10^{12}\);
\(\eta_C = 10^{-16}\);
\(j_{\max} = 5\).

For the nonmonotone strategies described in \Cref{sec:nonmonotone-strategies}, we choose
\(q = 12\);
\(r = 10\);
\(\eta_0 = \epsilonm^{0.5} \approx 10^{-8}\);
\(\eta_3 = 0.9\).
The limited memory BFGS approximation of \Cref{sec:quasi-newton-approximations} uses a default memory of \(6\).
For the small-step test of \Cref{sec:small-search-directions}, we choose \(n_{\mathrm{small}} = 2\).

\section{Software}%
\label{sec:software}

Penelopt, whose name derives from \emph{penalty optimization}, is written in Julia \citep{bezanson-edelman-karpinski-shah-2017} and is available as open-source software under the Mozilla Public License~2.0 \citep{penelopt}.
The results of this paper were obtained with version~1.0.0 of Penelopt and Julia~1.12.

Penelopt is part of the JSO ecosystem \citep{jso-2019}, which provides a common interface for nonlinear optimization solvers in Julia.
We use the NLPModels.jl API \citep{nlpmodels-2023} to model the optimization problems arising in our algorithm, and the SolverCore.jl API~\citep{solvercore-2023} for the implementation of the optimization algorithm.

Our continuous integration solves a set of small test problems at each commit, on Linux, macOS and Windows.
The documentation \citep{penelopt-docs} lists all parameters of the algorithm along with their default values, given in \Cref{sec:algorithmic-hyperparameters}.
It also contains tutorials that show how to call Penelopt on JuMP \citep{dunning-huchette-lubin-2017, lubin-dowson-diasgarcia-huchette-legat-vielma-2023}, AMPL \citep{fourer-gay-kernighan-1990}, and CUTEst \citep{gould-orban-toint-2015} models, and how to switch between different linear solvers.

For example, to solve the sixth problem from the Hock--Schittkowski collection \citep{hock-schittkowski-1981}, formulated as
\[
  \min_{x_1, x_2} {(1 - x_1)}^2 \quad \text{subject to } 10(x_2 - x_1^2) = 0,
\]
we can model it using the automatic differentiation interface of NLPModels.jl implemented in the ADNLPModels.jl package \citep{adnlpmodels-2021}.

\begin{lstlisting}[language=Julia, style=code]
using ADNLPModels

f(x) = (1 - x[1])^2
c(x) = [10 * (x[2] - x[1]^2)]

lcon = [0.0]
ucon = [0.0]

x0 = [-1.2, 1.0]

nlp = ADNLPModel(f, x0, c, lcon, ucon)
\end{lstlisting}

Then, we call the solver with the default algorithmic parameters,
\begin{lstlisting}[language=Julia, style=code]
using Penelopt
stats = L2Penalty(nlp; print_level = 1)
\end{lstlisting}

which results in the output shown in \Cref{lst:hs6-output}.
\begin{lstlisting}[style=output,
  caption={Output of the solver on HS6.},
  label={lst:hs6-output}]
This is Penelopt.jl v1.0.0.
Running with linear solver MUMPS v5.9.1.
------------------------------------------------------------------------------------------------------
Iter   sIter  Objective       pfeas       dfeas       τ           ptol        dtol        ‖x‖         
------------------------------------------------------------------------------------------------------
0      0      +4.8400000e+00  4.40e+00    1.56e+00    1.00e+00    1.00e+00    1.56e-02    1.56e+00    
1      8      +1.4515287e-04  4.25e-02    9.71e-03    1.00e+00    4.25e-04    9.71e-05    1.39e+00    
2      2      +9.9022550e-11  3.91e-06    7.88e-06    1.00e+00    8.05e-08    7.88e-08    1.41e+00    
3      1      +2.3024446e-15  9.81e-10    3.85e-08    1.00e+00    8.05e-08    3.82e-08    1.41e+00    
4      1      +5.7937134e-21  2.22e-14    6.09e-11    1.00e+00    8.05e-08    3.82e-08    1.41e+00    
Number of Iterations: 4
Objective...........: +5.793713424604307e-21
Primal Feasibility..:  2.220446049250313e-14
Dual Feasibility....:  6.090412449872034e-11
EXIT: first_order.
\end{lstlisting}

The solver reports the version and linear solver in use, followed by the iteration log and the final status.
On each iteration, the solver reports in order: the outer iteration number, the number of iterations of \Cref{alg:exactpen-inner-loop}, the objective value, the residuals \(\|c(x_k)\|_{\infty}\) and \(\|\nabla f_k + J_k^T y_k\|_{\infty}\), the penalty parameter \(\tau_k\), the primal and dual tolerances \(\epsilon_P^k\), \(\epsilon_D^k\), and the norm of the current iterate \(\|x_k\|_2\).
\section{Numerical Results}%
\label{sec:numerical-results}

In this section, we present how our algorithm compares with other state-of-the-art solvers on a set of benchmark problems.
Additionally, we provide numerical experiments to show the effect of the two nonmonotone strategies described in \Cref{sec:nonmonotone-strategies} on the performance of our algorithm.

We evaluate the performance of our algorithm on the CUTEst \citep{gould-orban-toint-2015} benchmark collection.
We select all CUTEst instances with equality constraints only, at least one constraint, at most as many constraints as variables, and no bounds other than fixed variables, which are removed in a preprocessing step.
This results in a set of \(224\) problems.
The largest problem in terms of variables has \(n = 132{,}098\) and \(m = 65{,}025\).
In terms of constraints, the largest problem has \(n = 123{,}200\) and \(m = 123{,}200\).

We compare our algorithm with the open-source IPOPT \citep{wachter-biegler-2006} v3.14 solver, which reduces to an SQP method for~\eqref{eq:nlp}.
As additional points of comparison, we also run MadNLP \citep{shin-anitescu-pacaud-2024} v0.11 and Uno v2.9 \citep{vanaret-leyffer-2026,vanaret-montoison-2026}.

All solvers run on each problem with a time limit of \(5\) minutes.
All benchmarks were run on GitHub-hosted \texttt{windows-latest} runners, which are x64 virtual machines running Windows Server with \(4\) virtual CPUs and \(16\)~GB of memory.
Moreover, all solvers declare optimality when~\eqref{eq:overall-stopping-criterion} is fulfilled with \(\epsilon_P = \epsilon_D = 10^{-6}\).

We disable the acceptable-level termination criteria of IPOPT, MadNLP, and Uno.
For all three of these, we turn off scaling of the dual feasibility measure \citep[(5)]{wachter-biegler-2006} by setting \(s_{\max}\) to the maximum representable Float64 value.
Although an optimization algorithm might encounter numerical difficulties when the Jacobian is ill-conditioned near a solution of~\eqref{eq:nlp}, we do not consider altering the stopping criterion an appropriate strategy.
First, this makes it difficult to compare the performance of different solvers.
Second, it gives a false impression of the quality of the solution to~\eqref{eq:nlp} when the problem is difficult.

All other parameters are set to their default values, and we use the hyperparameters from \Cref{sec:algorithmic-hyperparameters} for Penelopt.
All solvers, including Penelopt, use the MUMPS \citep{amestoy-duff-lexcellent-2000,amestoy-duff-lexcellent-koster-2001} v5.9 sparse linear solver to solve the linear systems arising in the computations.
We run Uno with its default \texttt{auto} preset, which selects the IPOPT preset when \(n + m \geq 2{,}000\) or when the Jacobian and the Hessian of the Lagrangian have at least \(50{,}000\) nonzeros in total, and the FilterSQP preset, with BQPD as the QP solver, otherwise.

Two exceptions apply to the problem set described above.
First, when using L-BFGS approximations, the problem FLOSP2TL caused the solvers to crash, hence we do not include it in L-BFGS benchmark runs.
Second, MadNLP and Uno repeatedly crashed on the problems CYCLOOCF and CYCLOOCT\@.
We therefore did not run these two solvers on these two problems, and we count both problems as failures for MadNLP and Uno.

\subsection*{CUTEst Benchmark}

To present the numerical results, we use \citet{dolan-more-2002} performance profiles.
Given a set \(\mathcal{P}\) of problems and a set \(\mathcal{S}\) of solvers, we compute a performance metric \(t_{p,s}\) for each \(p \in \mathcal{P}\) and \(s \in \mathcal{S}\).
We then compute the ratio \(r_{p,s} = t_{p,s} / \min_{s' \in \mathcal{S}} t_{p,s'}\) and define the performance profile of solver \(s\) as
\[
  \pi_s(\omega) \coloneqq \frac{1}{n_p} \sum_{p \in \mathcal{P}} \mathbf{1}_{r_{p,s} \leq \omega},
\]
where \(n_p = |\mathcal{P}|\), \(\omega \geq 1\), and \(\mathbf{1}_{r_{p,s} \leq \omega} = 1\) if \(r_{p,s} \leq \omega\) and \(0\) otherwise.

For the time profiles, we measure the wall-clock time \(t\), in seconds, reported by each solver.
Timings are recorded with a resolution of one millisecond.
To account for this timer resolution, we use \(t_{p,s} \coloneqq \max(t, 10^{-3})\) as the performance metric for time, so that two runs that both complete in under one millisecond are considered tied.

To test the infeasibility detection mechanism of \Cref{sec:infeasibility-detection}, we trigger an infeasibility certification procedure whenever a solver declares infeasibility.
Suppose a solver declares a problem infeasible at a point \(\bar x\).
To certify local infeasibility, we wish to find a local minimum \(x^*\) of
\[
  \min \|c(x)\|_2^2,
\]
that satisfies \(\|c(x^*)\|_2 > 0\) and \(\|J(\bar x) (\bar x - x^*)\|_2 \leq \Delta\).
In practice, we solve the equivalent problem
\[
  \min \|r\|_2^2 \quad \st \ c(x) - r = 0,
\]
using IPOPT with all default parameters, a time limit of \(5\) minutes, and a tolerance of \(10^{-9}\), started from \(x_0 = \bar x\).
We choose \(\Delta = 10\) and verify that \(\|c(x^*)\|_2 > 0\) by checking that \(\|c(x_{\mathrm{IPOPT}})\|_2 > 10^{-3}\), where \(x_{\mathrm{IPOPT}}\) is the solution returned by IPOPT\@.
Then, we check that \(\|J(\bar x) (\bar x - x^*)\|_2 \leq \Delta\) by verifying that \(\|J(\bar x) (\bar x - x_{\mathrm{IPOPT}})\|_2 \leq \Delta\).
If IPOPT fails to solve the problem or either condition is not met, we consider the certification procedure a failure.

Although this certification procedure might fail because \(x_{\mathrm{IPOPT}}\) is an approximation of \(x^*\), we consider it an effective measure of the quality of the infeasibility detection mechanism.
When a solver declares a problem infeasible, we consider the run successful if the certification procedure confirms local infeasibility and unsuccessful otherwise.
If both solvers declare infeasibility, the certification procedure may confirm only one of them.
In this case, we consider only the solver with confirmed infeasibility successful.

\begin{figure}[ht]%
  \centering
  \perfprofile{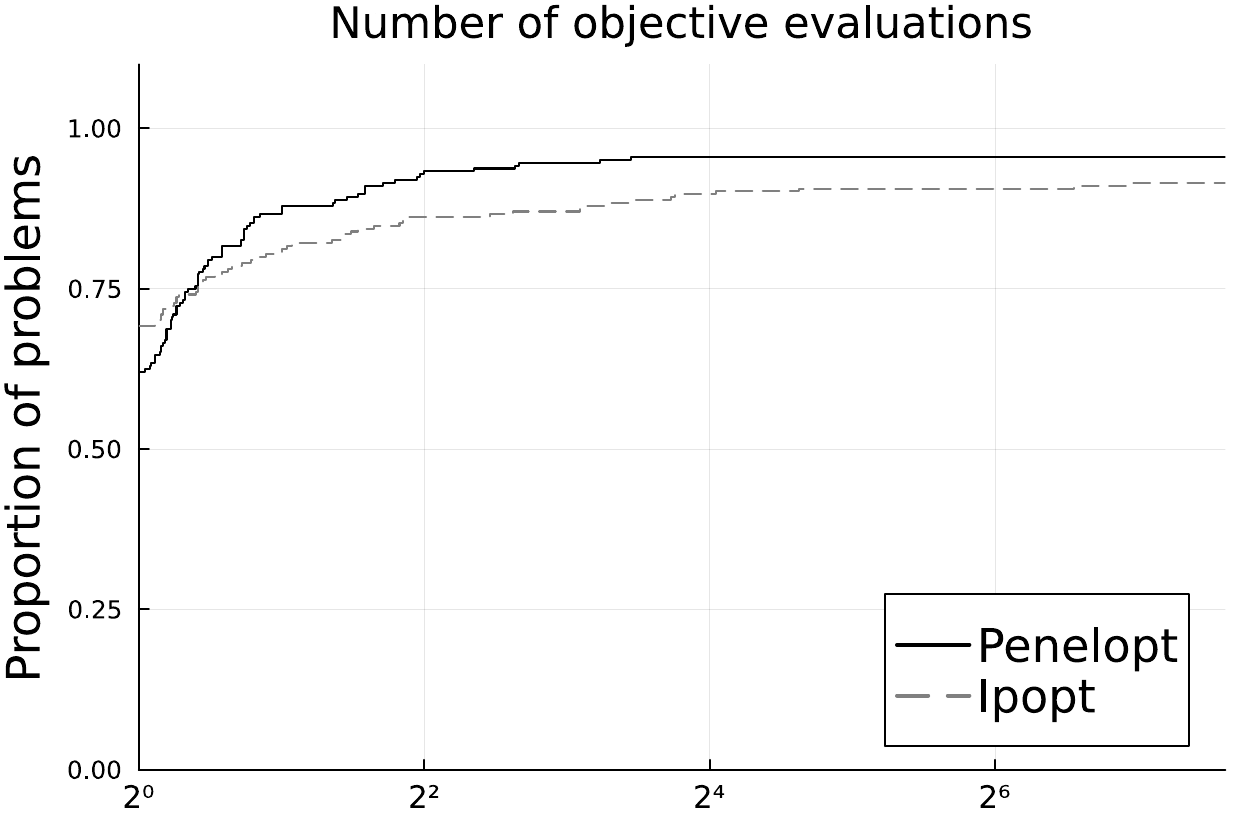}%
  \hspace{2em}%
  \perfprofile{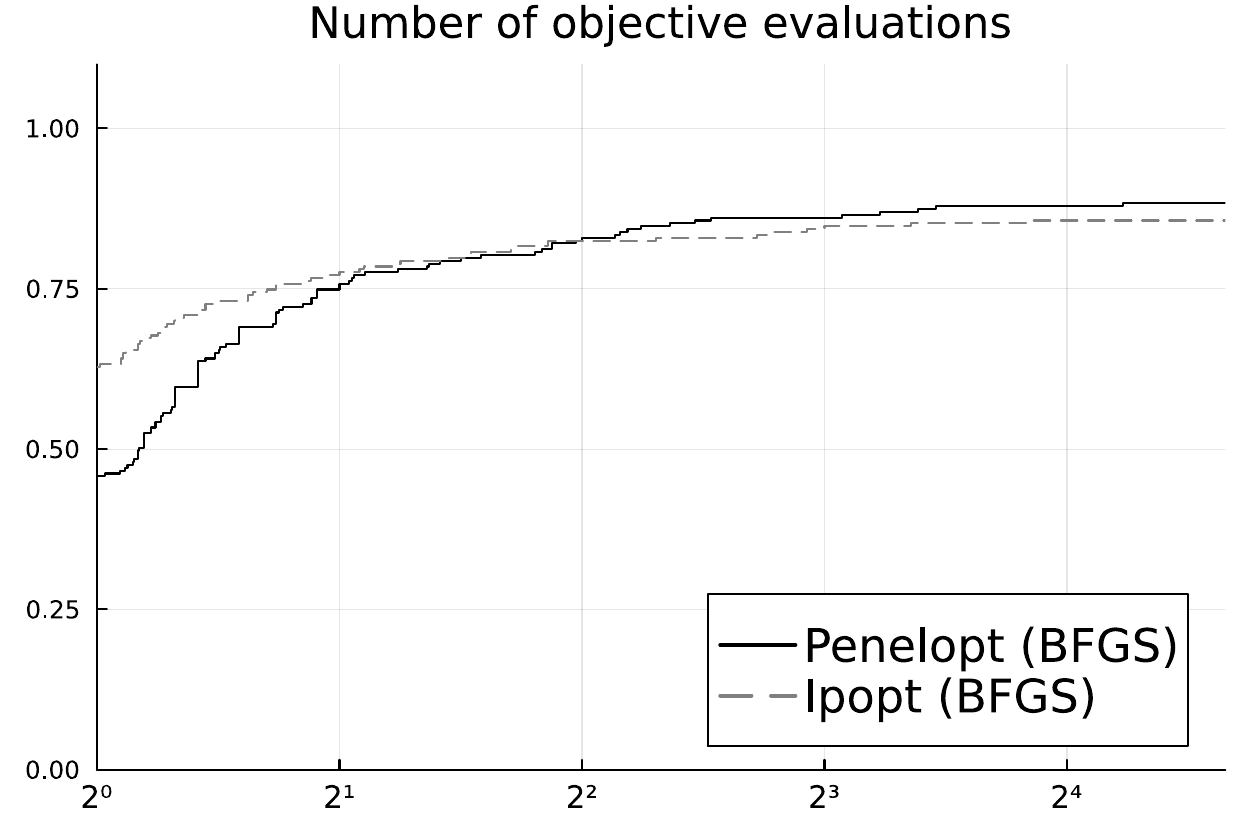}

  \medskip
  \perfprofile{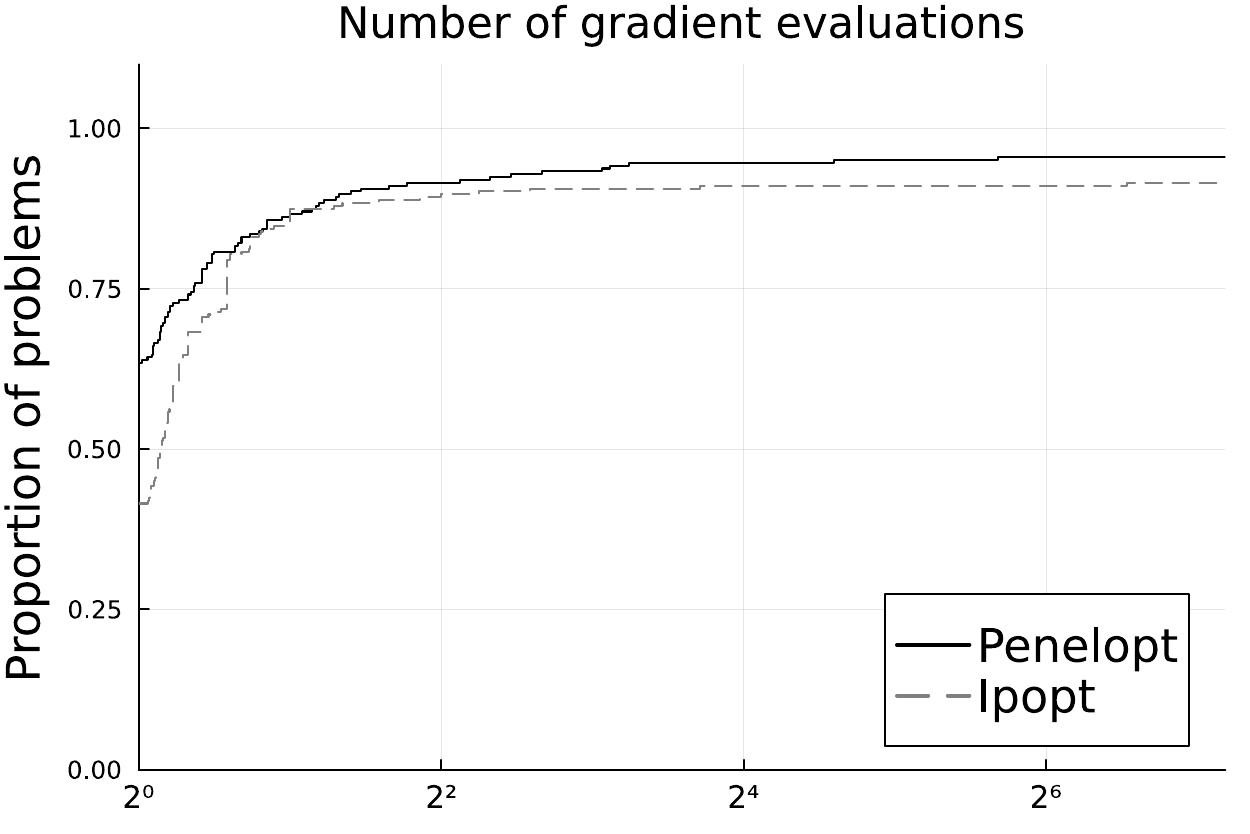}%
  \hspace{2em}%
  \perfprofile{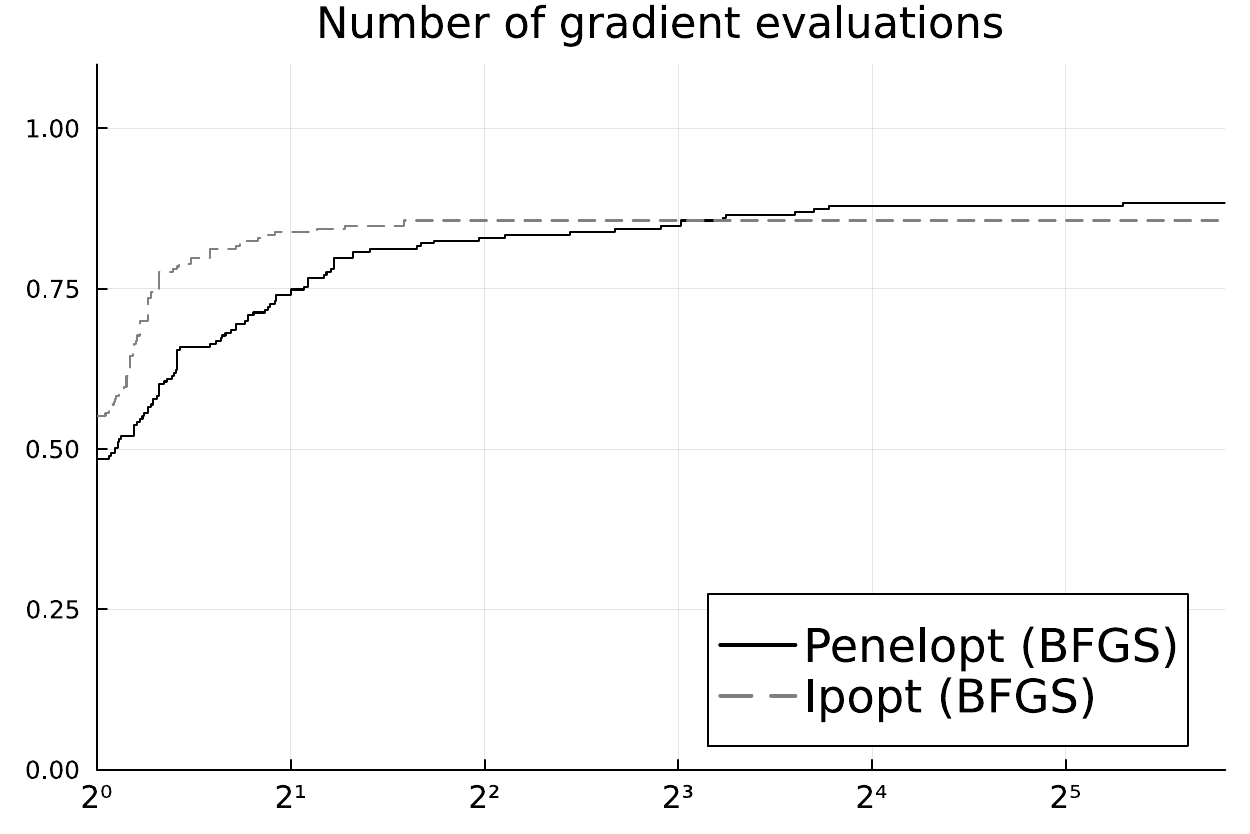}

  \medskip
  \perfprofile{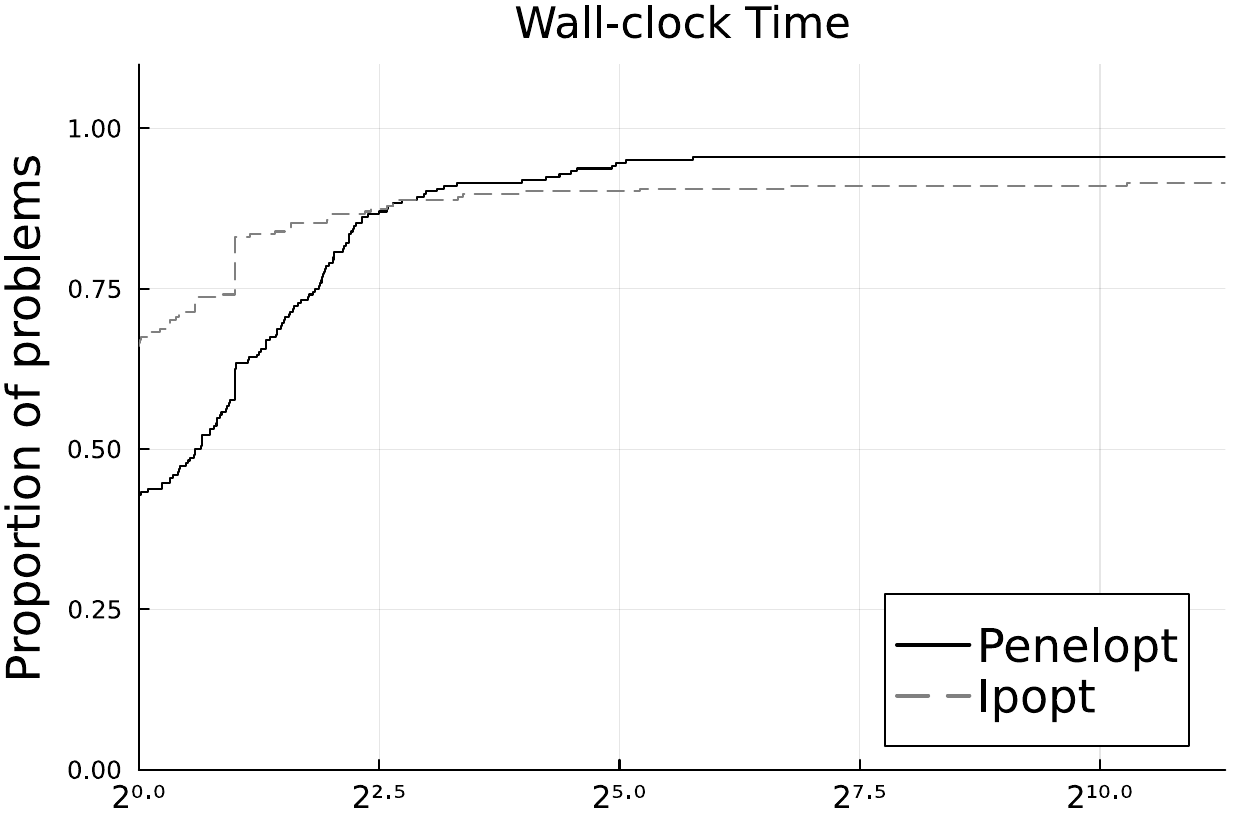}%
  \hspace{2em}%
  \perfprofile{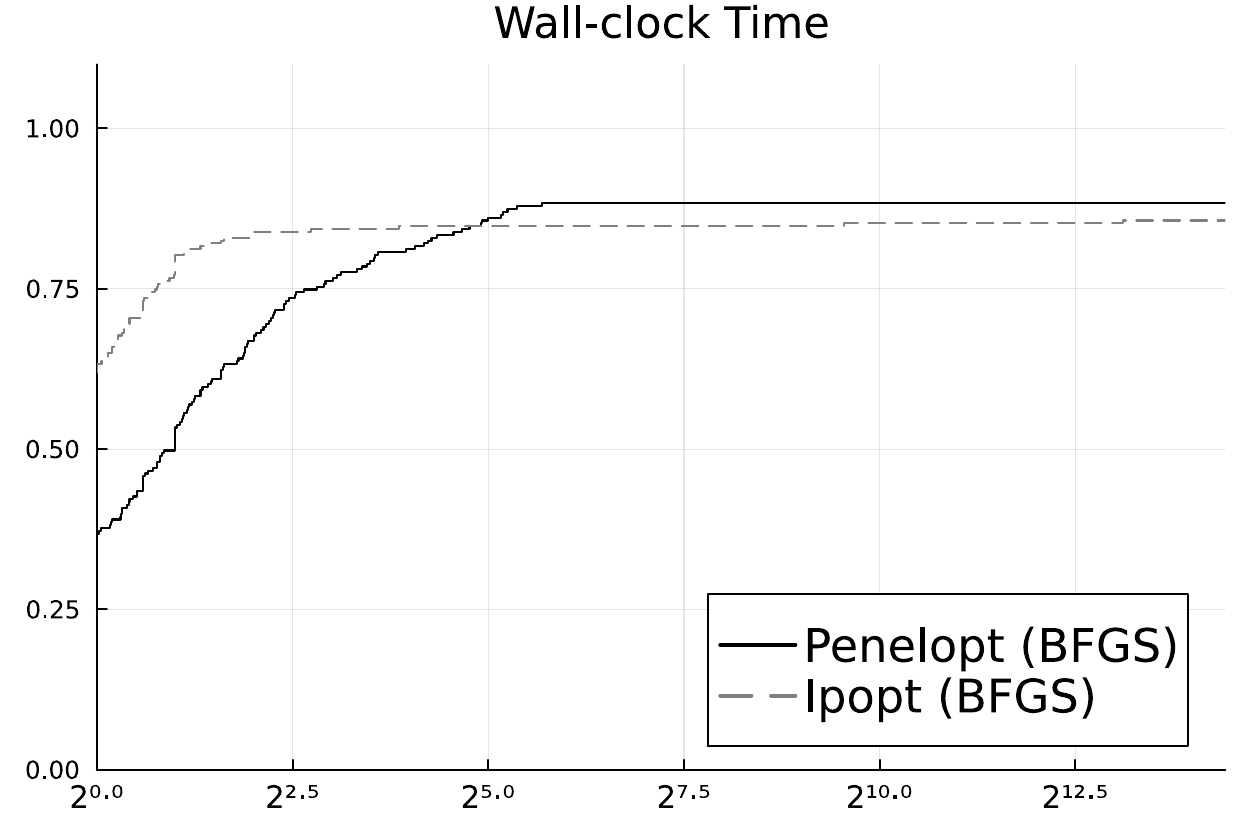}
  \caption{Left: Penelopt against IPOPT with \(H_l = \nabla^2 \mathcal{L}(x_l, y_l)\).
    Right: Penelopt against IPOPT with an L-BFGS quasi-Newton approximation of \(\nabla^2 \mathcal{L}(x_l, y_l)\), on \(223\) problems (FLOSP2TL excluded).
    Bottom row: wall-clock time, with times below \(1\)~ms set to \(1\)~ms.
    The higher, the better.}%
  \label{fig:perf-prof-ipopt}
\end{figure}

\begin{figure}[ht]%
  \centering
  \perfprofile{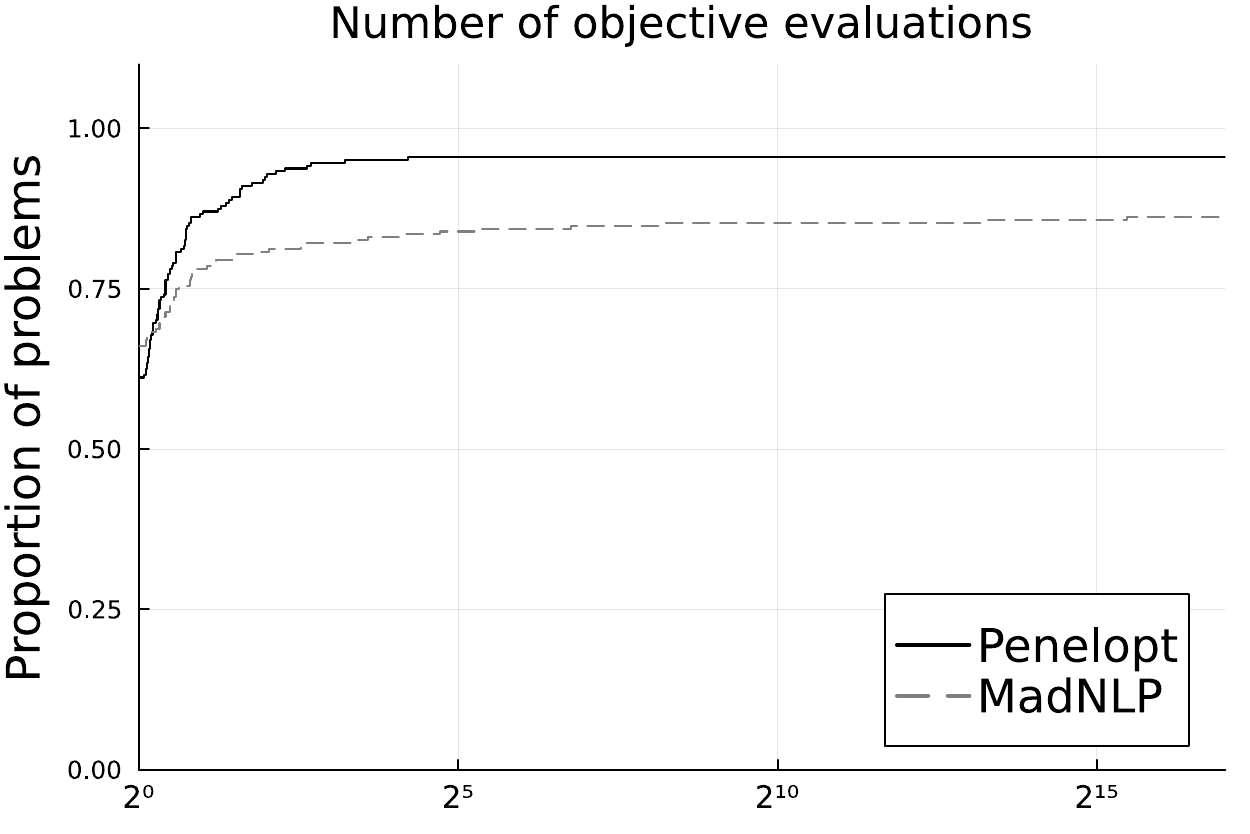}%
  \hspace{2em}%
  \perfprofile{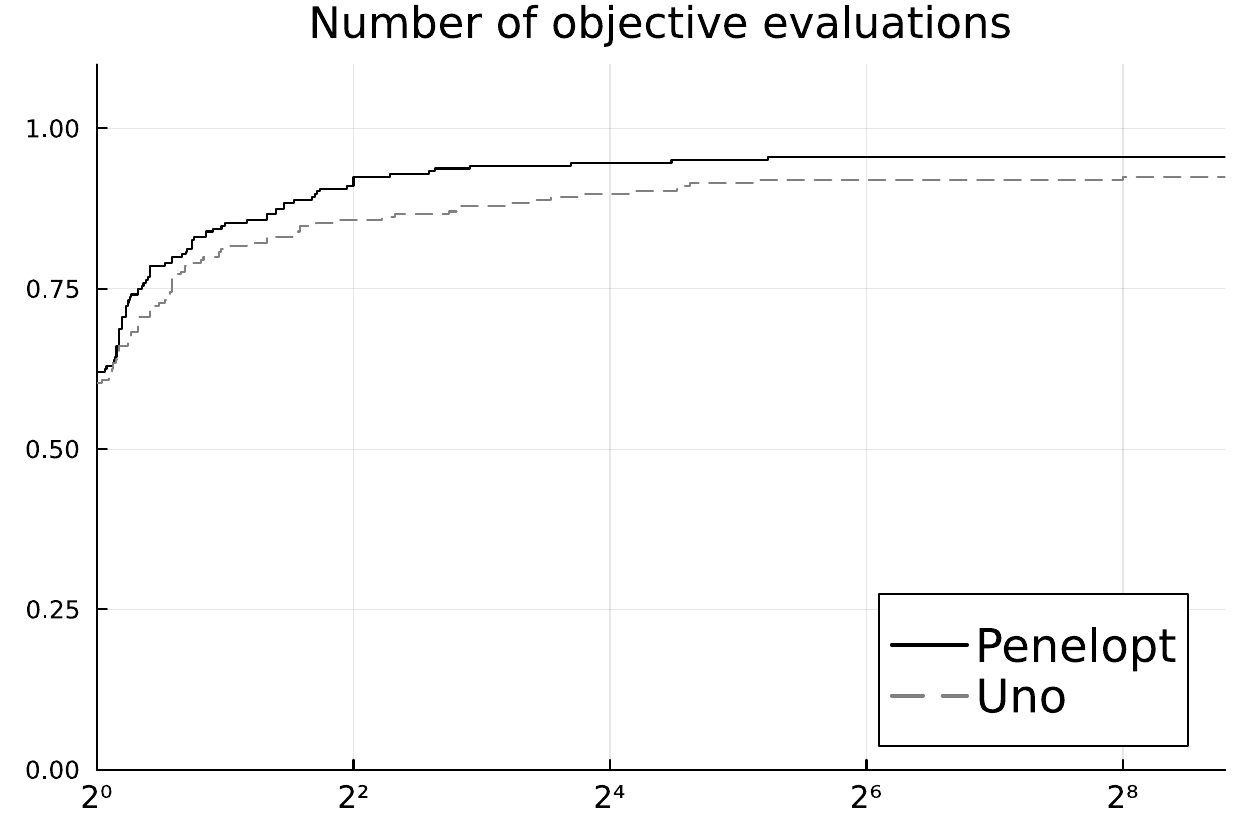}

  \medskip
  \perfprofile{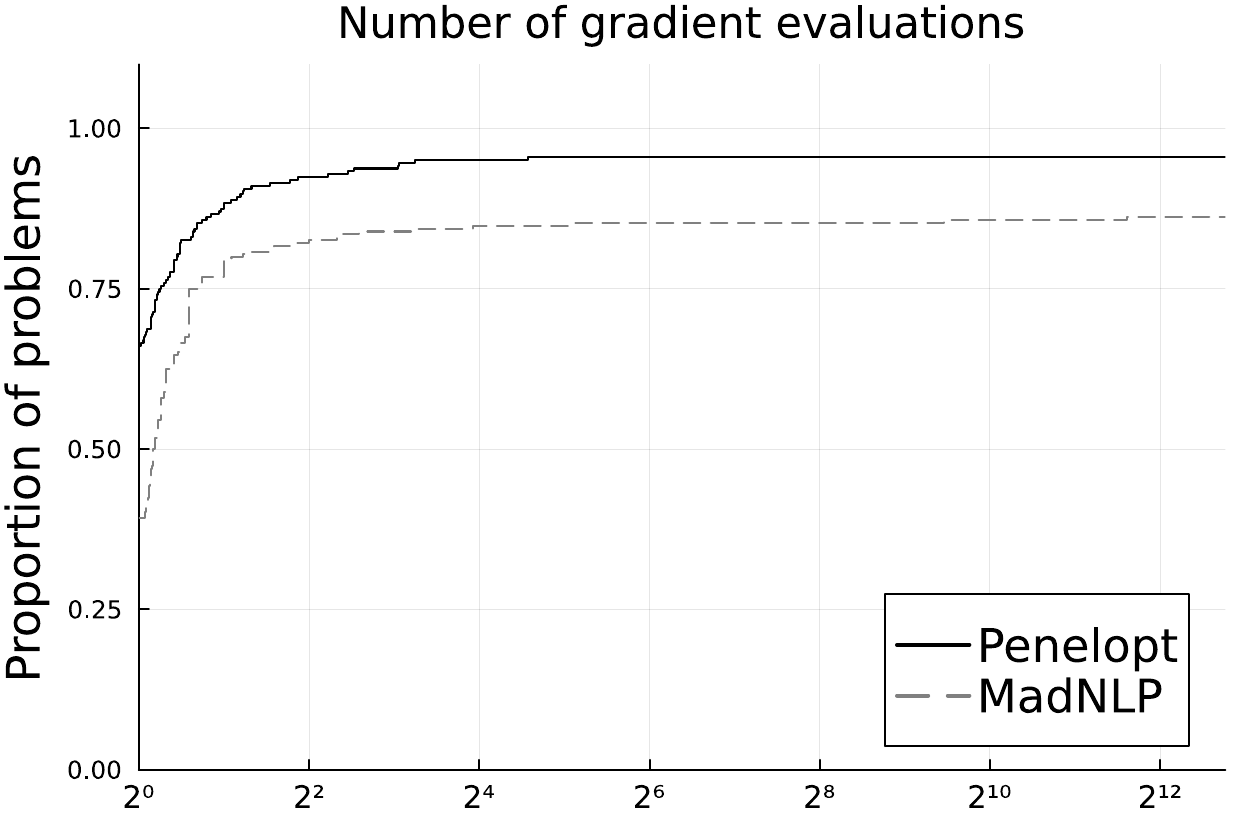}%
  \hspace{2em}%
  \perfprofile{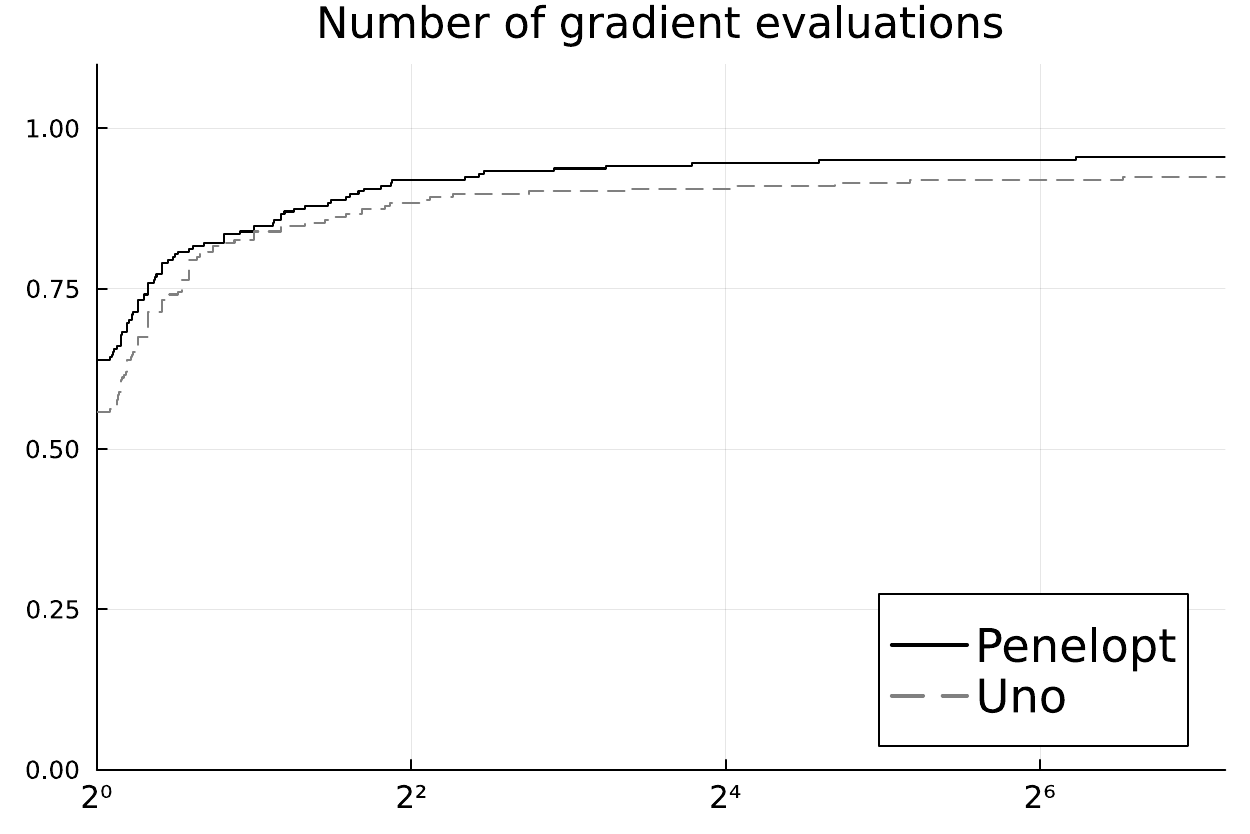}

  \medskip
  \perfprofile{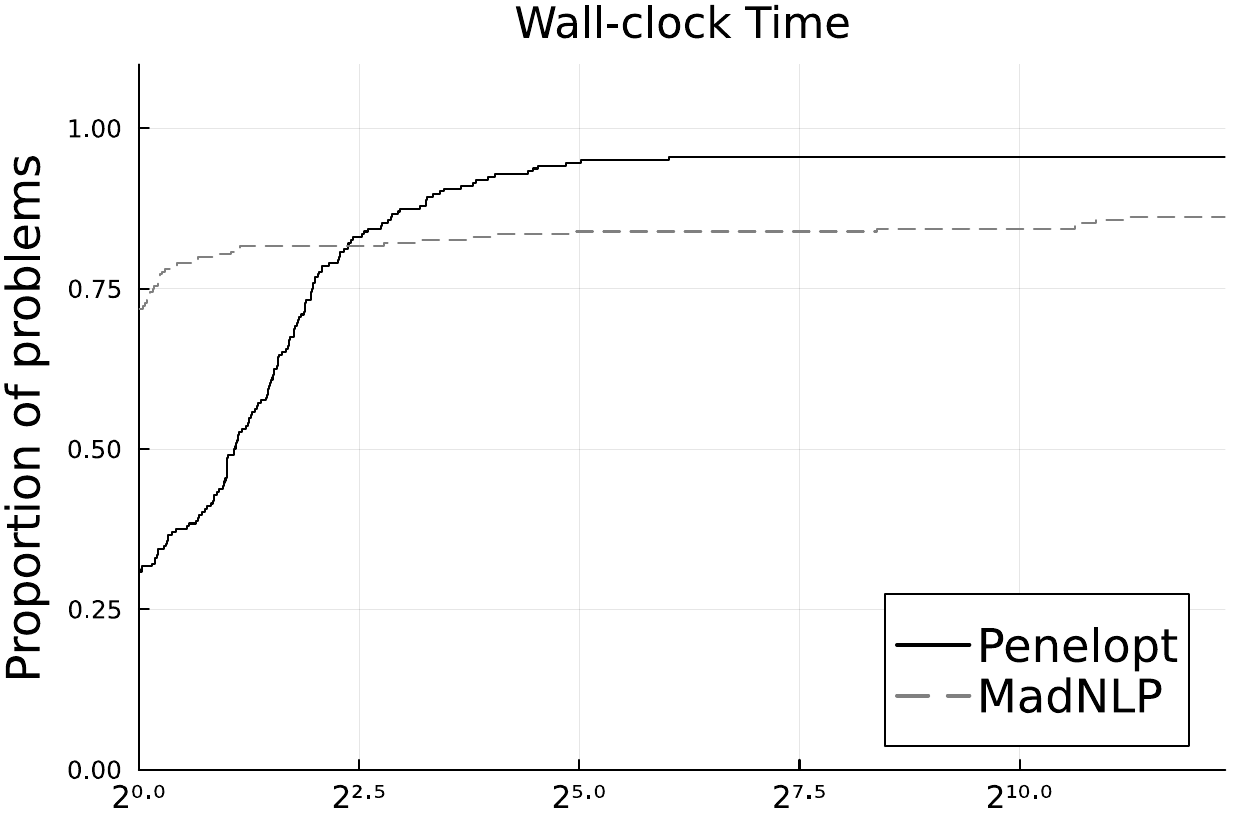}%
  \hspace{2em}%
  \perfprofile{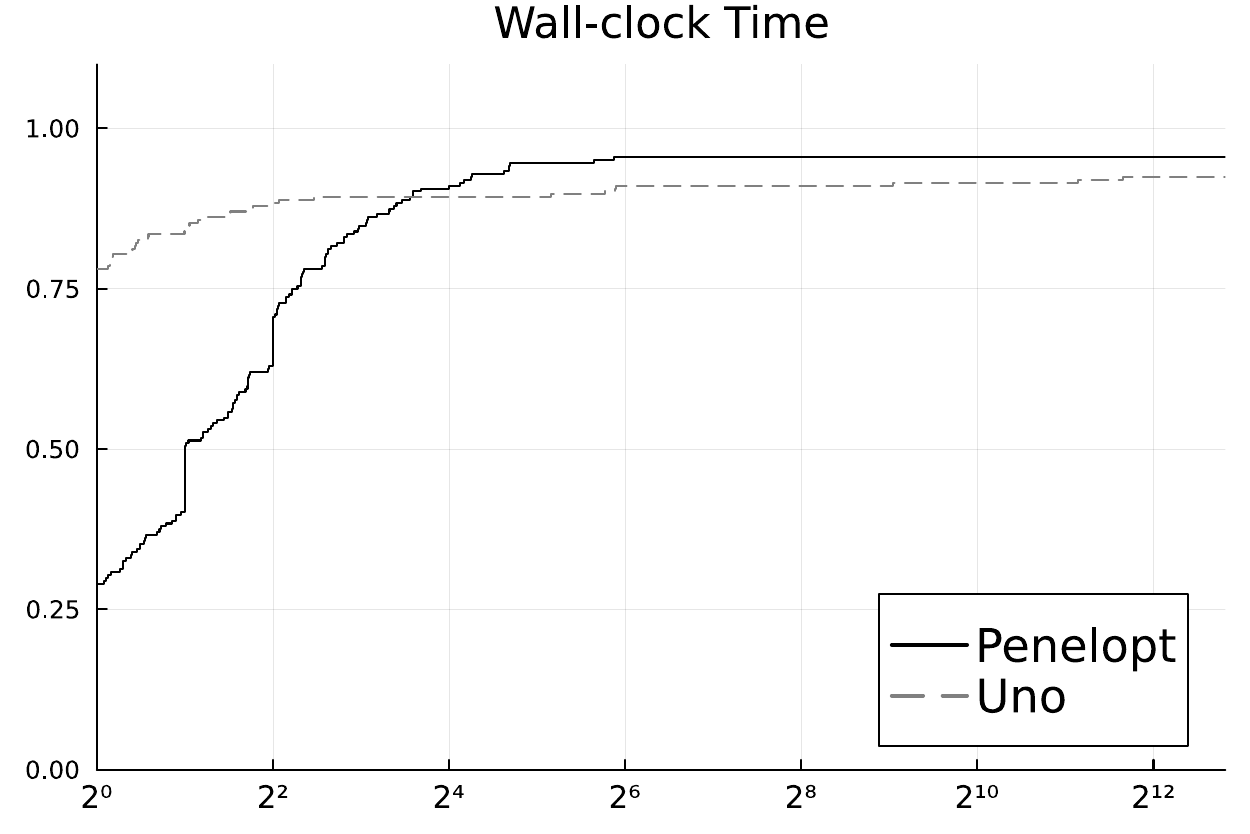}
  \caption{Left: Penelopt against MadNLP with \(H_l = \nabla^2 \mathcal{L}(x_l, y_l)\).
    Right: Penelopt against Uno with \(H_l = \nabla^2 \mathcal{L}(x_l, y_l)\).
    CYCLOOCF and CYCLOOCT count as failures for MadNLP and Uno.
    Bottom row: wall-clock time, with times below \(1\)~ms set to \(1\)~ms.
    The higher, the better.}%
  \label{fig:perf-prof-madnlp-uno}
\end{figure}

The source code to reproduce the numerical results is available in \citet{benchmark-runners}.
\Cref{fig:perf-prof-ipopt} reports performance profiles for Penelopt against IPOPT, with the exact Hessian and with an L-BFGS approximation.
\Cref{fig:perf-prof-madnlp-uno} additionally compares Penelopt, with \(H_l = \nabla^2 \mathcal{L}(x_l, y_l)\), against MadNLP and Uno under matching tolerances and time limits.

In \Cref{fig:perf-prof-ipopt}, we observe that Penelopt is competitive with IPOPT in terms of robustness and efficiency, with respect to the number of objective and gradient evaluations.
We see in \Cref{fig:perf-prof-ipopt} that Penelopt is slower than IPOPT in terms of wall-clock time.
We attribute this mainly to our step computation procedure (\Cref{alg:more-sorensen-step-indefinite}), which can require multiple factorizations of the KKT system~\eqref{eq:kkt-system} to compute a step \(s_l\).
Then, in \Cref{fig:perf-prof-madnlp-uno}, we observe that Penelopt is ahead of MadNLP and Uno in terms of robustness.

When \(H_l = \nabla^2 \mathcal{L}(x_l, y_l)\), Penelopt timed out on the problems EIGENBCO and DRCAVTY3, while it ended with a small step status on the problem CATENARY\@.
All other problems were either solved or declared infeasible.
In \Cref{tab:infeasibility-detection-exact}, we report the results of the infeasibility detection procedure described above.

\begin{table}[ht]%
  \centering
  \begin{tabular}{rcccccc}
  \hline
  \textbf{Problem} & \textbf{n} & \textbf{m} & \textbf{Penelopt (default)} & \textbf{Certified} & \textbf{IPOPT (default)} & \textbf{Certified} \\\hline
  ARTIF & 5002 & 5000 & infeasible & true & infeasible & true \\
  DENSCHNENE & 3 & 3 & first\_order &  & infeasible & true \\
  DRCAVTY3 & 4489 & 3969 & max\_time &  & infeasible & true \\
  FLOSP2HH & 2883 & 2763 & infeasible & false & infeasible & false \\
  FLOSP2HL & 2883 & 2763 & infeasible & true & infeasible & true \\
  FLOSP2HM & 2883 & 2763 & infeasible & true & infeasible & true \\
  FLOSP2TH & 2883 & 2763 & infeasible & false & first\_order &  \\
  FREURONE & 2 & 2 & infeasible & true & infeasible & true \\
  HATFLDFLNE & 3 & 3 & infeasible & false & first\_order &  \\
  HIMMELBD & 2 & 2 & infeasible & true & infeasible & true \\
  JUNKTURN & 10010 & 7000 & first\_order &  & infeasible & false \\
  LINCONT & 1257 & 419 & infeasible & true & infeasible & false \\
  LUKVLE15 & 9997 & 7497 & first\_order &  & infeasible & false \\
  LUKVLE16 & 9997 & 7497 & first\_order &  & infeasible & false \\
  NONDIANE & 5000 & 5000 & infeasible & true & max\_time &  \\
  NONMSQRTNE & 4900 & 4900 & infeasible & false & infeasible & false \\
  OSCIPANE & 10 & 10 & infeasible & false & first\_order &  \\
  POROUS1 & 5184 & 4900 & infeasible & false & first\_order &  \\
  POWELLSQ & 2 & 2 & infeasible & true & infeasible & true \\
  POWERSUMNE & 4 & 4 & first\_order &  & infeasible & false \\
  VANDANIUMS & 22 & 10 & infeasible & true & unknown &  \\
  WOODSNE & 4000 & 3001 & infeasible & true & infeasible & true \\
  YATP2CNE & 123200 & 123200 & infeasible & false & first\_order &  \\\hline
\end{tabular}
  \caption{Infeasibility detection results on the CUTEst benchmark with \(H_l \coloneqq \nabla^2 \mathcal{L}(x_l, y_l)\).
    Columns: problem name, number of variables, number of constraints, Penelopt status, infeasibility certification for the Penelopt solution, IPOPT status, infeasibility certification for the IPOPT solution.
    False indicates that we were unable to certify infeasibility.}%
  \label{tab:infeasibility-detection-exact}
\end{table}

\subsection*{Effect of Nonmonotone Strategies}

\begin{figure}[!tbp]%
  \centering
  \perfprofile{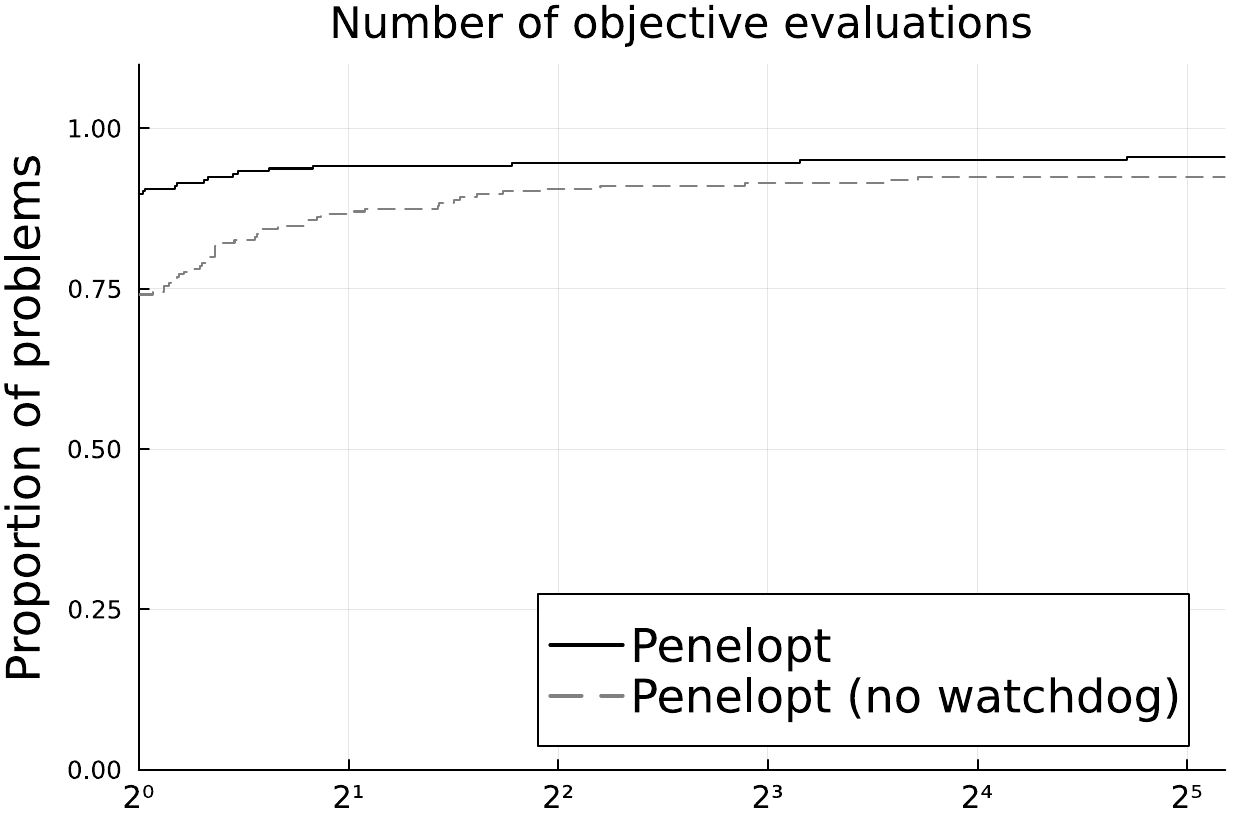}%
  \hspace{2em}%
  \perfprofile{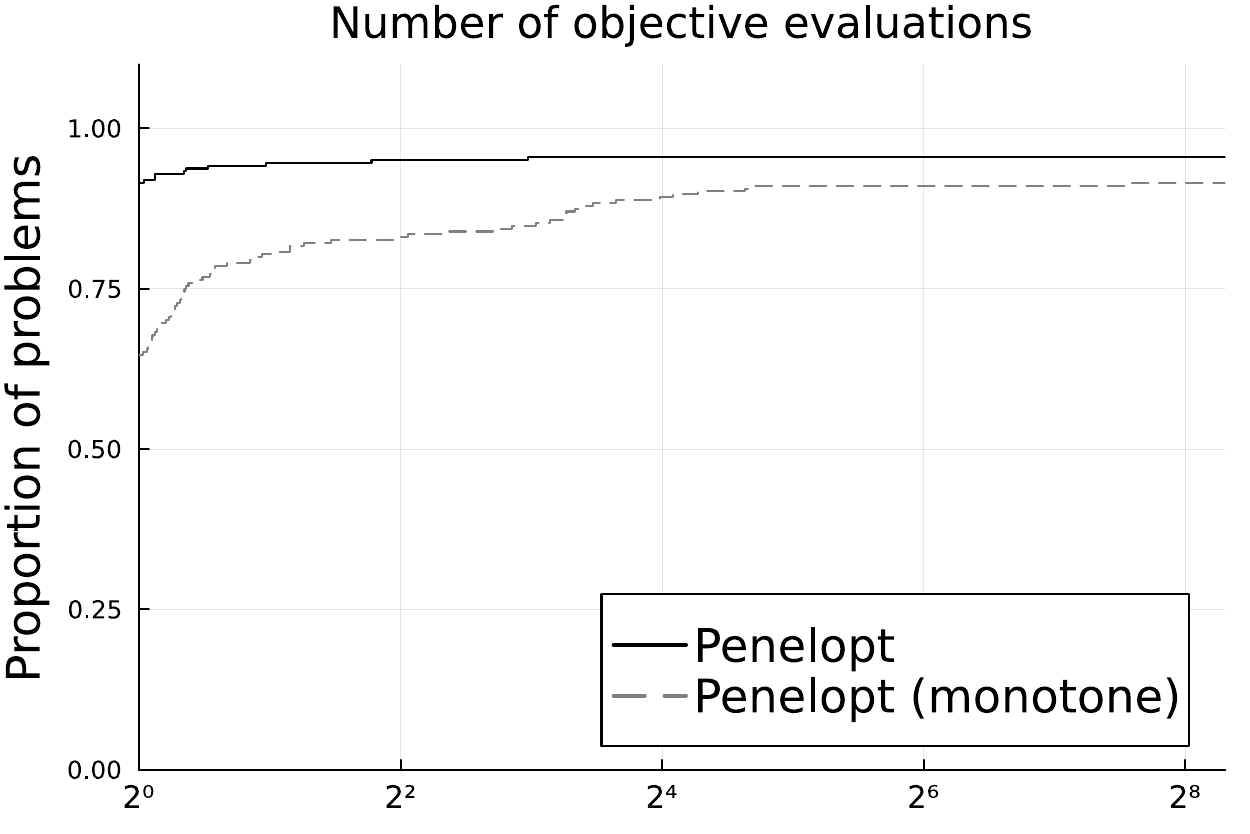}

  \medskip
  \perfprofile{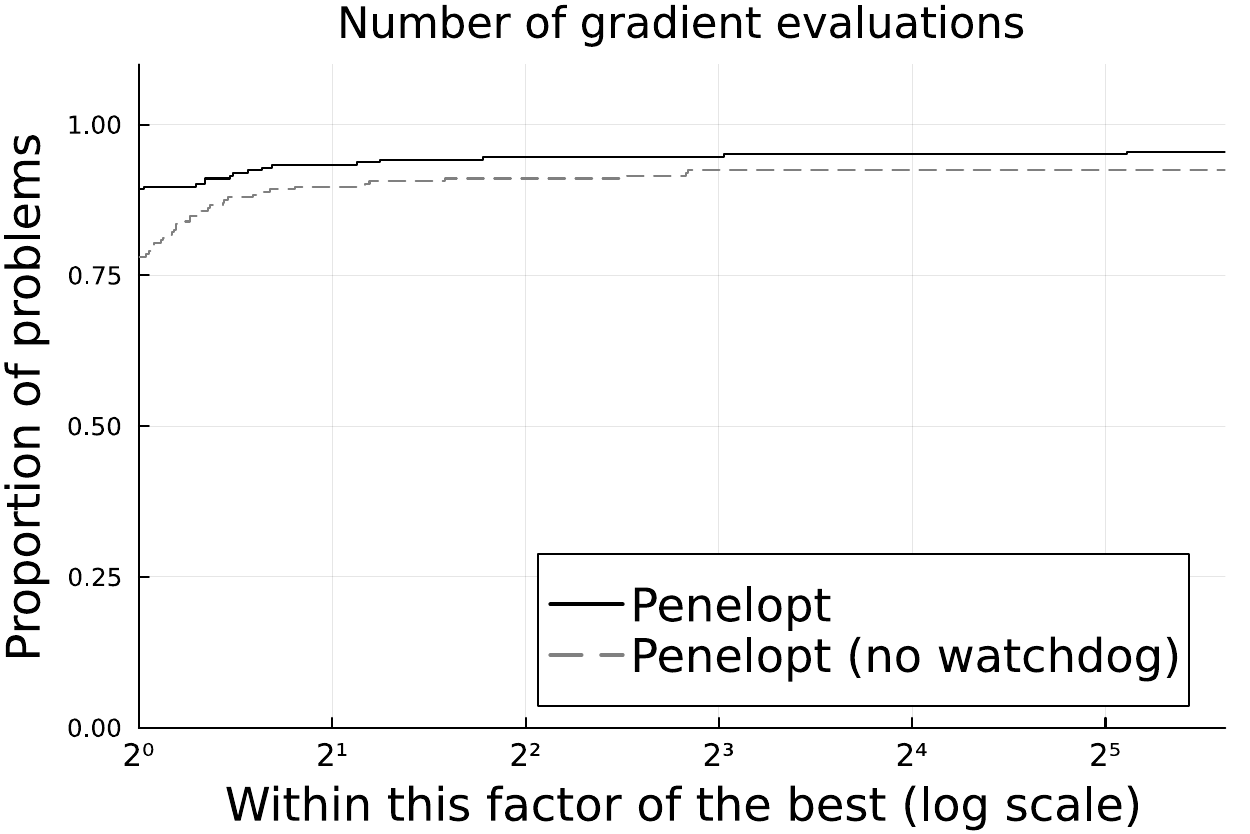}%
  \hspace{2em}%
  \perfprofile{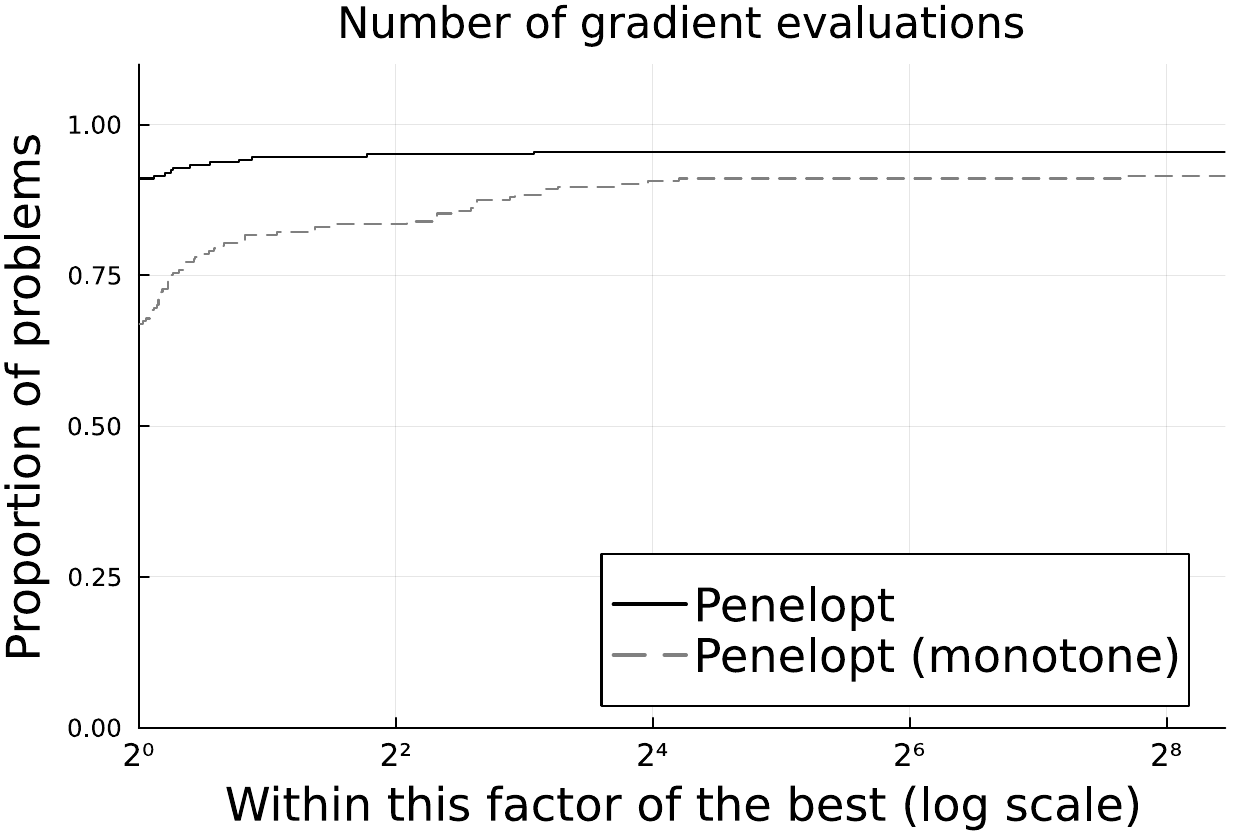}
  \caption{Ablation study with \(H_l = \nabla^2 \mathcal{L}(x_l, y_l)\).
    Left: Penelopt against Penelopt without the watchdog.
    Right: Penelopt against Penelopt with the monotone acceptance rule~\eqref{eq:armijo}.
    The higher, the better.}%
  \label{fig:ablation}
\end{figure}

We now present the effect of the two strategies of \Cref{sec:nonmonotone-strategies}, which, among the heuristics of \Cref{sec:implementation-details}, have the largest impact on the performance of Penelopt in our experiments.
We compare Penelopt, with \(H_l = \nabla^2 \mathcal{L}(x_l, y_l)\) and the default hyperparameters of \Cref{sec:algorithmic-hyperparameters}, with two variants: one in which the watchdog is disabled, and one in which steps are accepted with the monotone ratio~\eqref{eq:armijo} instead of~\eqref{eq:armijo-nonmonotone}, the watchdog being kept.
All other settings are as above.
\Cref{fig:ablation} reports the corresponding performance profiles in terms of objective and gradient evaluations.

We see on \Cref{fig:ablation} that both nonmonotone strategies improve the robustness of Penelopt, with the update rule~\eqref{eq:armijo-nonmonotone} having a larger impact than the watchdog rule.

\section{Conclusion}%
\label{sec:conclusion}

Penelopt is the first large-scale implementation of the exact penalty method.
On large-scale problems, our implementation is competitive with IPOPT and places the exact penalty method on par with SQP methods in terms of robustness and efficiency.

The extension of the method to handle bound constraints is the subject of active research.
In particular, we are investigating the use of a logarithmic barrier function to handle inequalities.

We are also considering alternative linear algebra routines for solving the linear systems arising in \Cref{alg:more-sorensen-step-definite,alg:more-sorensen-step-indefinite}.
Specifically, we are exploring preconditioned iterative methods to replace the sparse \(LDL^T\) factorizations.
In addition, we plan to evaluate Penelopt in single and quadruple precision.
For this purpose, in addition to MUMPS and MA57, our solver is compatible with LDLFactorizations.jl \citep{ldlfactorizations-2023}, which is also part of the JSO ecosystem, implements a sparse \(LDL^T\) factorization in pure Julia and supports single and quadruple precision.
In a similar fashion, we plan to extend our implementation to GPU architectures by using cuDSS \citep{cudss,montoison-2024} in place of MUMPS and MA57 for the linear algebra routines, as in MadNLP\@.

The superlinear local convergence rate of \Cref{alg:exactpen-inner-loop} also remains to be established.
Additionally, it remains to prove that \Cref{alg:exactpen-outer-loop} converges to a solution of~\eqref{eq:nlp} when the penalty parameter \(\tau_k\) is unbounded.
Finally, a more theoretically grounded strategy for the indefinite \(H_l + \sigma_l I\) case of \Cref{alg:more-sorensen-step-indefinite} is also a subject of ongoing research.

\section*{AI Usage Disclosure}%
Generative AI tools were used to improve the language and grammar of this paper, to screen the text for mathematical inconsistencies, and to draft and revise parts of the text; all suggested references were verified by the authors.
Claude (Anthropic) was used to assist with code refactoring and documentation during the development of Penelopt.
The authors reviewed and edited all such output as needed and take full responsibility for the content of this paper.

\bibliographystyle{abbrvnat}
\bibliography{abbrv, penalty-implementation}
\end{document}